\documentclass[11pt]{amsart}
\usepackage{amsmath}
\usepackage{amsthm}
\usepackage{amsfonts}
\usepackage{amssymb}
\usepackage[all]{xy}
\usepackage{alltt}
\usepackage{enumitem}
\usepackage{stmaryrd}
\usepackage{comment}
\usepackage{xspace}
\usepackage{comment}
\usepackage[
 colorlinks=true,urlcolor=blue,linkcolor=blue,citecolor=blue
]{hyperref}

\theoremstyle{plain}
\newtheorem{prop}{Proposition}[section]

\newtheorem{problem}{Problem}[section]

\newtheorem{corollary}[prop]{Corollary}

\newtheorem{thm}[prop]{Theorem}
\theoremstyle{definition}

\newtheorem{lem}[prop]{Lemma}
\newtheorem{remark}[prop]{Remark}
\newtheorem{example}[prop]{Example}

\newcommand{\Z}{\ensuremath{\mathbf{Z}}}
\newcommand{\F}{\ensuremath{\mathbf{F}}}

\newcommand{\Q}{\ensuremath{\mathbf{Q}}}

\newcommand{\C}{\ensuremath{\mathbf{C}}}
\newcommand{\GL}{\ensuremath{\mathbf{GL}}}
\newcommand{\SL}{\ensuremath{\mathbf{SL}}}
\newcommand{\PGL}{\ensuremath{\mathbf{PGL}}}
\newcommand{\PSL}{\ensuremath{\mathbf{PSL}}}

\newcommand{\n}{\noindent}

\renewcommand{\F}{\ensuremath{\mathbf{F}}}
\renewcommand{\Z}{\ensuremath{\mathbf{Z}}}
\renewcommand{\C}{\ensuremath{\mathbf{C}}}
\newcommand{\A}{\ensuremath{\mathbf{A}}}
\newcommand{\D}{\ensuremath{\mathbf{D}}}

\usepackage[margin=1.25in]{geometry}

\usepackage[
        colorlinks=true,urlcolor=blue,linkcolor=blue,citecolor=blue
]{hyperref}

\begin{document}  
 
\title{Fuchs' problem for endomorphisms \\ of nonabelian groups}
\date{\today}
 
\author{Sunil K. Chebolu}
\address{Department of Mathematics \\
Illinois State University \\
Normal, IL 61790, USA}
\email{schebol@ilstu.edu}

\author{Keir Lockridge} 
\address {Department of Mathematics \\
Gettysburg College \\
Gettysburg, PA 17325, USA}
\email{klockrid@gettysburg.edu}

\thanks{Sunil Chebolu is partially supported by Simons Foundation's Collaboration Grant for Mathematicians (516354)}

\keywords{group of units, group algebra, quaternion group, dihedral group, symmetric group, fully realizable, Fuchs' problem}
\subjclass[2020]{Primary 16U60, 08A35; Secondary 16S34}
 
\begin{abstract}

In 1960, L\'{a}szl\'{o} Fuchs posed the problem of determining which groups $G$ are realizable as the group of units in some ring $R$. In \cite{chebolu2022fuchs}, we investigated the following variant of Fuchs' problem, for abelian groups: which groups $G$ are realized by a ring $R$ where every group endomorphism of $G$ is induced by a ring endomorphism of $R$? Such groups are called fully realizable. In this paper, we answer the aforementioned question for several families of nonabelian groups: symmetric, dihedral, quaternion, alternating, and simple groups; almost cyclic $p$-groups; and groups whose Sylow $2$-subgroup is either cyclic or normal and abelian. We construct three infinite families of fully realizable nonabelian groups using iterated semidirect products.
\end{abstract}
 
\maketitle
\thispagestyle{empty}

\tableofcontents

\section{Introduction}
Given a ring $R$ with identity, let $R^\times$ denote the group of units in $R$, and for any ring endomorphism $\phi$ of $R$, let $\phi^\times$ denote the induced group endomorphism of $R^\times$. This paper is a sequel to \cite{chebolu2022fuchs} in which we introduced the following problem: 

\begin{problem}Find all groups $G$ for which there is a ring $R$ such that\label{prob}
\begin{itemize}
    \item[(1)] $R^\times = G$, and
    \item[(2)] every group endomorphism of $G$ is of the form $\phi^{\times}$ for some  ring endomorphism $\phi$ of $R$.
\end{itemize}
\end{problem}

\n If condition (1) is satisfied, then we say the ring $R$ {\bf realizes} $G$ and $G$ is {\bf realizable}. If both conditions are satisfied, then we say the ring $R$ {\bf fully realizes} $G$ and the group $G$ is {\bf fully realizable}. Problem \ref{prob} is a generalization of one originally posed by L\'{a}szl\'{o} Fuchs: which groups $G$ are realizable as the group of units in some ring $R$? Given that the assignment $R \longmapsto R^\times$ from rings to groups is functorial, it seems natural to also consider maps and ask when both the group {\em and} all of its endomorphisms are realizable. We refer the reader to the introduction of \cite{chebolu2022fuchs} for more on the genesis of and history behind Fuchs' problem and this variant.

In \cite{chebolu2022fuchs}, we solved Problem \ref{prob} for groups of odd order, torsion-free abelian groups, torsion abelian groups, and finitely generated abelian groups. (If a group of odd order is the group of units in a ring $R$, then the subring of $R$ generated by its units is commutative, so \cite{chebolu2022fuchs} was almost exclusively a study of full realizability in the category of abelian groups.) Our objective in this paper is to address Problem \ref{prob} for the class of nonabelian groups. In particular, we show that there are no nonabelian fully realizable groups of the following types: dihedral groups, split extensions of a cyclic group by a cyclic group of prime order, quaternion groups, finite simple groups, symmetric groups, alternating groups (on five or more symbols), groups with a cyclic Sylow 2-subgroup, and finite general linear groups. In \cite{chebolu2019fuchs}, we proved that the only realizable nonabelian almost cyclic $p$-groups (this class includes the dihedral, quaternion, semidihedral, and modular 2-groups) are the dihedral group of order 8 and the quaternion group of order 8, neither of which is fully realizable, so there are no fully realizable nonabelian almost cyclic $p$-groups. 

Fortunately, we are able to give three infinite families of positive examples. Our examples are of the form $A \rtimes \C_{p^n}$ where $p^n$ is a prime power and $A$ is a finite 2-group. Theorem \ref{ea2gpn} summarizes our positive cases when $A$ is an abelian 2-group, and Theorem \ref{c6} (together with the surrounding discussion) provides a family of positive examples where $A$ is a nonabelian 2-group. Concretely, we construct the following families of fully realizable nonabelian groups as iterated semidirect products; $\C_n$ denotes the cyclic group of order $n$ and $\A_4$ is the alternating group on 4 symbols.
\begin{itemize}
    \item Let $A = \C_2 \times \C_2$, $G_1 = A \rtimes \C_4$ (GAP id (16, 3)) and \[G_i = \underbrace{A \rtimes  A \rtimes \cdots \rtimes A}_i\rtimes \C_4\] for $i \geq 2$.
    \item Let $A = \C_2 \times \C_2$, $G_1 = A \rtimes \C_3 = \A_4$ and \[G_i = \underbrace{A \rtimes  A \rtimes \cdots \rtimes A}_i\rtimes \C_3\] for $i \geq 2$.
    \item Let $A = \C_2^4 \rtimes \C_2$ (GAP id (32, 27)), $G_1 = A \rtimes \C_3$ (GAP id (96, 70)), and \[G_i = \underbrace{A \rtimes A \rtimes \cdots \rtimes A}_i\rtimes \C_3\] for $i \geq 2$.
\end{itemize}

\n The case $p^n = 3$ is special, and we are able to classify the fully realizable groups $G$ where the Sylow 2-subgroup of $G$ is normal and abelian (see Theorem \ref{na2sg}), because in this case $G = A \rtimes \C_3$, with $A$ an abelian 2-group, is forced (see Proposition \ref{normSyl2}).

We use GAP (\cite{GAP4}) and SageMath (\cite{sagemath}) to prove that certain specific finite groups either are or are not fully realizable. Our code is contained in several Jupyter notebook files in the directory \cite[\href{https://cocalc.com/klockrid/main/FPENAG}{FPENAG}]{ourcode}; the auxiliary functions we used are contained in the file \cite[\href{https://cocalc.com/klockrid/main/FPENAG/files/FPENAG.py}{FPENAG.py}]{ourcode}. We also found the GroupNames website to be an indispensable guide to the classification of groups of small order and their properties; see \cite{GN}.

Semidirect products $A\rtimes B$ are used extensively in this paper. If the homomorphism $\phi \colon B \longrightarrow \mathrm{Aut}(A)$ that determines the semidirect product is not clear from context, then our intent is that $A\rtimes B$ stand for an arbitrary semidirect product of $A$ and $B$ (this is often the case in the statements of our theorems involving semidirect products). We will write $\phi_b$ for $\phi(b)$ when convenient.


\noindent {\bf Acknowledgement}. We are
 grateful to the reviewer for providing detailed and helpful comments.
\section{Preliminaries}

We begin by collecting a few results from our earlier work and giving several definitions. In \cite[3.1]{chebolu2022fuchs}, we proved that if a ring $R$ fully realizes its group of units $G$, then $R$ has characteristic 2 and $G$ is in fact fully realized by the subring $S$ of $R$ that is generated by $G$. Further, this subring $S$ is isomorphic to $\F_2[G]/I$ for some (two-sided) ideal $I$ in the group ring $\F_2[G]$.

Call an ideal $I$ in a group ring $k[G]$ {\bf $\boldsymbol{\mathrm{End}(G)}$-invariant} if, for every endomorphism $\phi$ of $G$, the natural extension $\overline{\phi}\colon k[G] \longrightarrow k[G]$ maps $I$ into $I$. If the restriction of the quotient map $k[G] \longrightarrow k[G]/I$ to $G$ is injective, we say that $G$ {\bf embeds} in $k[G]/I$. In \cite[3.1 and 3.2]{chebolu2022fuchs} we proved that if there is an ideal $I \subseteq \F_2[G]$ such that (1) $G$ embeds in $\F_2[G]/I$, (2) $I$ is $\mathrm{End}(G)$-invariant, and (3) $\F_2[G]/I$ realizes $G$, then the ring $\F_2[G]/I$ fully realizes $G$. We will therefore say that an ideal $I$ satisfying these three conditions fully realizes $G$. 

A subgroup $H$ of a group $G$ is a retract if there is a homomorphism from $G$ to $H$ whose restriction to $H$ is the identity map. The kernel of this homomorphism is a normal complement of $H$ in $G$, and in fact $H$ is a retract of $G$ if and only if $H$ has a normal complement. If $S$ is a set of ring elements, then $(S)$ will denote the two-sided ideal generated by the elements in $S$; if $S$ is a set of group elements, then $\langle S \rangle$ will denote the subgroup generated by the elements in $S$.

\begin{prop} \label{prprop} \phantom{E}
\begin{enumerate}
    \item The fully realizable groups of odd order are $\C_1$ and $\C_3$ {\em(\cite[4.1]{chebolu2022fuchs})}.
    \item A finite cyclic group $\C_n$ is fully realizable if and only if $n$ is a divisor of $12$ {\em (\cite[6.4]{chebolu2022fuchs})}. \label{realizablecyclics}
    \item The fully realizable finite abelian groups are the groups $W \times H$ where $W$ is a finite elementary abelian $2$-group and $H$ is a subgroup of $\C_{12}$ {\em (\cite[6.1]{chebolu2022fuchs})}. \label{realizablefab}
    \item A group $G$ is fully realizable if there is an ideal $I \subseteq \F_2[G]$ such that $I$ fully realizes $G$. Any ring $R$ generated by its units that fully realizes $G$ is isomorphic to $\F_2[G]/I$ for some ideal $I$ that fully realizes $G$ {\em (\cite[3.1, 3.2]{chebolu2022fuchs})}. \label{I-invariant}
    \item Every retract of a fully realizable group is also fully realizable; thus, if a semidirect product  $A \rtimes B$ of groups is fully realizable, then so is $B$. In particular, direct summands of fully realizable groups are fully realizable {\em (\cite[3.3]{chebolu2022fuchs})}. \label{PQsd} 
    \item If $G$ and $H$ are fully realizable groups and there are no nontrivial homomorphisms from $G$ to $H$ or from $H$ to $G$, then $G \times H$ is fully realizable {\em (\cite[3.4]{chebolu2022fuchs})}. 
\end{enumerate}
\end{prop}

\begin{example}
Let $\mathbf{V}_4 = \C_2 \times \C_2 = \langle a \rangle \times \langle b \rangle$ denote the Klein-four group. One can check that $\mathbf{V}_4$ embeds in $R = \F_2[\mathbf{V}_4]/(1 + a + b + ab)$ and $R^\times$ is isomorphic to $\mathbf{V}_4$. Moreover, the ideal $I = (1 + a + b + ab)$ is  an $\mathrm{End}(\mathbf{V}_4)$-invariant ideal. To see this, let $f$ be any endomorphism of $\mathbf{V}_4$. If $f$ is trivial, then $\bar{f}(1+a+b+ab) = 1+1+1+1 = 0$. If $f$ is an automorphism, then it permutes the 4 group elements, so $\bar{f}(1+a+b+ab) = 1+a+b+ab$. Finally, if the kernel of $f$ is a proper, nontrivial subgroup of $\mathbf{V}_4$, we may assume without loss of generality (by changing basis) that the kernel is $\langle a \rangle$. Then, $\bar{f}(1+a+b+ab) = 1 + 1 + f(b) +f(b) =0$. This shows that $I$ is $\mathrm{End}(\mathbf{V}_4)$-invariant and hence $\mathbf{V}_4$ is fully realizable by Proposition \ref{prprop} (\ref{I-invariant}). This argument can be extended to arbitrary elementary abelian $2$-groups, so all such groups are fully realizable.
\end{example}

\begin{remark} As in the example above, elements of the form $1 + x + y + xy$ appear frequently in this paper. This is essentially a consequence of \cite[3.8]{chebolu2022fuchs}. In short, if $S$ is the subring of the group ring $\F_2[G] \times \F_2[H]$ generated by $G \times H$, then $S$ is isomorphic to the quotient ring $\F_2[G \times H]/I,$ where $I$ is the ideal generated by all elements of the form $(1 + x)(1+y) = 1 + x + y + xy,$ where $x \in G$ and $y \in H.$\end{remark}

\begin{example} \label{dihedralex}
Let $n = 2k$ be even with $k \geq 3$ odd. The dihedral group \[\D_{2n} = \langle r,  s \,|\, r^n =1, s^2 =1, srs^{-1} = r^{-1} \rangle\] has direct product decomposition $\D_{2n} \cong \C_2 \times \D_{n}.$ Hence, in this case, if $\D_{2n}$ is fully realizable, then $\D_{n}$ is fully realizable by Proposition \ref{prprop} (\ref{PQsd}). (If $k$ is even, then $\D_{2k}$ is not a retract of $\D_{4k},$ so Proposition \ref{prprop} (\ref{PQsd}) does not apply.)

\end{example}

\begin{example} \label{generallinear}
For $n\ge 1$ and $q$ a prime power, the determinant map on $\GL_n(\F_q)$ gives the following split short exact sequence of groups:

\[
\xymatrix{
1 \ar[r] & \SL_n(\F_q) \ar[r] & \GL_n(\F_q)   \ar[r]^{\ \ \text{det}}   &  \F_q^{\times} \ar@/^1pc/[l]^{\omega}\ar[r] & 1,
}
\]
where $\omega(x)$ is the $n \times n$ diagonal matrix whose $(i, i)$ entry is $x$ when $i=1$ and $1$ when $i > 1$. Thus, $\GL_n(\F_q) = \SL_n(\F_q) \rtimes (\F_q)^\times$, and hence $(\F_q)^\times \cong \C_{q-1}$ is a retract of $\GL_n(\F_q)$. Therefore, if $\GL_n(\F_q)$ is fully realizable, then so is 
$\C_{q-1}$. By Proposition \ref{prprop} (\ref{realizablecyclics}), $q-1$ must divide $12$, so $q$ is $2 ,3, 4, 5, 7$ or $13$. In Section \ref{sec:fglg}, we use this example to find the fully realizable finite general linear groups.
\end{example}

The next three propositions will be used to study certain split extensions of cyclic groups (see \S \ref{sec:cbyc}).

\begin{prop} \label{maxabelian}
Let $G$ be a group and let $H$ be a maximal abelian subgroup of $G$ with the property that every endomorphism of $H$ can be extended to an endomorphism of $G$. If $G$ is fully realizable, then so is $H$.  
\end{prop}

\begin{proof}
Let $G$ and $H$ be as given in the statement of the proposition and let $G$ be fully realized by a ring $R$, necessarily of characteristic $2$. Consider the subring $S$ of $R$ generated by $H$. Note that $H \subseteq S^\times$ and  $S$ is commutative. The maximality of $H$ therefore implies that $S^\times = H$. Let $\phi \colon H \rightarrow H $ be any endomorphism. Then there exists an endomorphism $\rho \colon G \rightarrow G$ such that $\rho|_H = \phi$.  Since $G$ is fully realizable, there exists a map $\bar{\phi} \colon R \longrightarrow R$ such that $\bar{\phi}^\times = \rho$. Furthermore, $\bar{\phi}|_S$ sends $S$ to $S$ and    $(\bar{\phi}|_S)^\times = \phi$.
\end{proof}

\begin{example}
If $\D_{2n}$ is fully realizable, then $\C_n$ is fully realizable. It is easy to see that any endomorphism $f \colon \C_n \longrightarrow \C_n$ that sends $r \mapsto r^d$ (where $d$ is a divisor of $n$) can be extended to an endomorphism of $\D_{2n}$ by sending $s$ to itself. Therefore, if $\D_{2n}$ is fully realizable, then so is $\C_n$, and hence  $n$ divides $12$. In Section \ref{sec:dihedral} we will prove that there are no fully relizable nonabelian dihedral groups; our proof will not rely on the classification of fully realizable cyclic groups.
\end{example}

\begin{prop}
Let $G = A \rtimes B$ be a semidirect product of groups corresponding to a homomorphism $\phi \colon B \rightarrow Aut(A)$. An endomorphism $f \colon A \longrightarrow A$ can be extended to an endomorphism $\hat{f} \colon G \rightarrow G$ if $f h = h f$ for all $h$ in the image of $\phi$.
\end{prop}

\begin{proof}
Let $f$ be an endomorphism of $A$, and assume that $f h = h f$ for all $h$ in the image of $\phi$. Define $\hat{f} \colon G \rightarrow G$ by sending $(a, b) \mapsto (f(a), b)$. Then we have the following.

\begin{eqnarray*}
\hat{f}( (a_1, b_1)(a_2, b_2)) & = & \hat{f}(a_1 \phi_{b_1}(a_2), b_1b_2) \\
& = & (f(a_1 \phi_{b_1}(a_2)), b_1b_2) \\
& = & (f(a_1)f(\phi_{b_1}(a_2)), b_1b_2)\\
& = & (f(a_1)\phi_{b_1}(f(a_2)), b_1b_2) \\
& = & (f(a_1), b_1)(f(a_2), b_2)\\
& = & \hat{f}(a_1, b_1) \hat{f}(a_2, b_2)
\end{eqnarray*}
This shows that $\hat{f}$ is an endomorphism of $G$ that extends $f$. 
\end{proof}

\begin{prop} \label{max-ab-semi}
Let $A$ be an abelian group and let $G = A \rtimes B$ be a semidirect product of groups corresponding to a homomorphism $\phi \colon B \rightarrow Aut(A)$. Further assume that $\text{Im}(\phi) \subseteq Z(\text{End}(A))$ and $A$ is a maximal abelian subgroup of $G$ or $\phi$ is trivial. If $G$ is fully realizable, then $A$ is fully realizable. 
\end{prop}

\begin{proof}
If $\phi$ is trivial, then $A$ is fully realizable by Proposition \ref{prprop} (\ref{PQsd}). Otherwise, the given condition $\text{Im}(\phi) \subseteq Z(\text{End}(A))$ guarantees that every endomorphism  on $A$ can be extended to an endomorphism of $G$, and then we can apply Proposition \ref{maxabelian} to complete the proof.
\end{proof}

Finally, we prove a proposition that will be used to study two cases: (1) groups with a cyclic Sylow 2-subgroup and (2) groups whose Sylow 2-subgroup is normal. Let $G$ and $T$ be groups. Any homomorphism $\phi \colon G \longrightarrow \mathrm{Aut}(T)$ induces an action of $G$ on $T$ defined by $t^g = \phi(g)(t)$ for $t \in T$ and $g\in G$. An element $t \in T$ is fixed by $G$ if $t^g = t$ for all $g \in G$. The set $T_f$ of elements in $T$ that are fixed by $G$ is a subgroup of $T$. If $T_f = \{1\}$, then we say the action is {\bf fixed-point free}. If $G$ acts on a product of groups $A \times B$, then we say the action is {\bf decomposable} (and $G$ {\bf acts decomposably}) if $a^g \in A$ and $b^g \in B$ for all $a \in A, b \in B,$ and $g \in G$.

\begin{prop} \label{fpfsplit}
Let $p$ be a prime, let $H$ be a finite abelian group with $\gcd(|H|, p) = 1$, and suppose $\C_{p^k} = \langle x \rangle$ acts on $H$ via an automorphism $h \mapsto h^x$ of order $p$. Then, there is a decomposition $H = H_f \times H'$ where $\C_{p^k}$ acts decomposably on $H_f \times H'$ and the action of $\C_{p^k}$ on $H'$ is fixed-point free.
\end{prop}
\begin{proof}
Since $H$ is a finite abelian group, we may write it as a direct sum of $q$-groups $H_q$ where each $q \neq p$ is prime: $H = \oplus_q H_q$. Since $x$ acts on $H$ via an automorphism, it preserves the orders of elements, so $x$ acts decomposably on this direct sum. Fix a summand $H_q$ and let $q^k$ be the exponent of $H_q$. Then, $H_q $ is a $\Z_{q^k}[\C_p] = \Z_{q^k}[x]/(x^p - 1)$-module since $x$ acts via an automorphism of order $p$.

Note that $x^p - 1 = (x - 1)\Phi_p(x)$ where \[\Phi_p(x) = 1 + x + \cdots + x^{p-1}.\] Since $x - 1$ is monic, we may divide $\Phi_p(x)$ by $x - 1$ to obtain \[\Phi_p(x) = h(x)(x-1) + \alpha,\] where $\alpha = \Phi_p(1) = p$ is a unit in $\Z_{q^k}$. Hence, \[1 = (-1/\alpha)h(x)(x-1) + (1/\alpha)\Phi_p(x).\] This proves that $(x-1, \Phi_p(x)) = \Z_{q^k}[x].$ Further, $(x -1) \cap (\Phi_p(x)) = (x^p - 1)$ because $1$ is not a root of $\Phi_p(x)$. By the Chinese Remainder Theorem for rings, \[\Z_{q^k}[\C_p] = \Z_{q^k} \times \Z_{q^k}[x]/(\Phi_p(x)).\] We may now write $H_q = A_q \times B_q,$ where $x$ acts trivially on the $\Z_{q^k}$-module $A_q$ and $B_q$ is a $\Z_{q^k}[x]/(\Phi_p(x))$-module. The action of $x$ on $A_q$ is trivial, and if $xm = m$ for some $m \in B_q$ then \[0 = \Phi_p(x)m = pm \implies m = 0\] since $p$ is a unit mod $q^k$. Thus, the action of $x$ on $B_q$ is fixed-point free. We now have $H_f = \oplus_q A_q$ and, defining $H' = \oplus_q B_q$, the action of $\C_{p^k}$ on $H'$ is fixed-point free. This completes the proof.
\end{proof}

The conclusion of Proposition \ref{fpfsplit} need not hold if $p \mid |H|$ (e.g., for the dihedral group $\C_4 \rtimes \C_2$ of order 8, $H = \C_4$, $H_f = \C_2$, and $H_f$ is not a summand). The conclusion also need not hold if $H$ is not abelian; for example, let $H = \C_2^4 \rtimes \C_2$ (GAP id (32, 27)) and let $G = H \rtimes \C_3$ (GAP id (96, 70)). In this example, $H_f = \C_2$ is not a summand of $H$. We will later discuss the fact that $G$ is fully realizable but $H$ is not; see \cite[\href{https://cocalc.com/klockrid/main/FPENAG/files/SG-96-70.ipynb}{SG-96-70.ipynb}]{ourcode}.

\section{Dihedral groups}\label{sec:dihedral}

The realizable dihedral groups were determined in \cite{chebolu2017fuchs}. In positive characteristic, only the groups $\D_k$ for $k = 2,4,6, 8, 12$ are realizable, and in characteristic zero, only $\D_2$ and $\D_{4k}$ for $k$ odd are realizable. Throughout this section, we will use the notation established in Example \ref{dihedralex}. Our goal is to prove the following theorem.

\begin{thm} \label{dihedral}
The dihedral group $\D_{2n}$ is fully realizable if and only if $n = 1$ or $2$.
\end{thm}

\n We separately consider the cases where $n$ is odd, $n = 2k$ for $k$ odd, and $4 \mid n$.

\begin{prop}
If $4 \mid n$, then $\D_{2n}$ is not fully realizable.
\end{prop}
\begin{proof}
Assume to the contrary that there is a ring $R$ of characteristic 2, generated by its units, that fully realizes $\D_{2n}$. For any $l$, there is an endomorphism $\phi_l \in \mathrm{End}(\D_{2n})$ defined by $\phi_l(r) = r$ and $\phi_l(s) = sr^l$. 

Suppose $4 \mid n$. Then, \[
\begin{aligned}
(1 + s + sr^{n/4})^2 &= r^{n/4} + r^{-n/4} + 1\\
(1 + s + sr^{n/4})^4 &= 1 + r^{n/2} + r^{-n/2} = 1.
\end{aligned}
\] Since $1 + r^{n/4} + r^{-n/4} \neq 1$, $1 + s + sr^{n/4}$ must have order 4 as a unit of $R$. Thus, $1 + s + sr^{n/4} = r^k$ for some $k$. If we apply $\phi_{n/4}$ to this relation, we obtain that $1 + sr^{n/4} + sr^{n/2} = r^k$. Adding these equations yields $r^{n/2} = 1$, a contradiction.
\end{proof}

\begin{prop}
If $n \geq 3$ is odd, then $\D_{2n}$ is not fully realizable. If $n = 2k$ for $k\geq 3$ odd, then $\D_{2n}$ is not fully realizable.
\end{prop}
\begin{proof}
Let $n \geq 3$ be odd. Assume to the contrary that there is a ring $R$ of characteristic 2, generated by its units, that fully realizes $\D_{2n}$. Consider the element \[ \sigma = 1 + r + \cdots + r^{n-1}\] of $R$. Since $|r|$ is odd, $\langle r \rangle = \langle r^2 \rangle$ and thus $\sigma^2 = \sigma$ (the element $\sigma$ is idempotent). Since $\langle r \rangle = \langle r^{-1} \rangle$, $\sigma$ is central. Thus $\sigma$ and $\tau = 1 + \sigma$ are orthogonal central idempotents of $R$ and $R$ decomposes as a direct product of rings $R\sigma \times R\tau$. However, $\D_{2n}$ is indecomposable, so one of these factors has a trivial unit group.

If $(R\sigma)^\times = \{1\}$, then \[s\sigma = \sigma.\] Consider the endomorphism of $\D_{2n}$ defined by $\phi(r) = 1$ and $\phi(s) = s$. Applying this to the previous equation, we obtain $s = 1$, a contradiction.

If $(R\tau)^\times = \{1\}$, then $r\tau = \tau$. By direct computation, $r\sigma = \sigma$. From this it follows that $r = r(1) = r(\sigma + \tau) = r\sigma + r\tau = \sigma + \tau = 1$, a contradiction.

Now suppose $n = 2k$ where $k\geq 3$ is odd. In this case, $\D_{4k} \cong \C_2 \times \D_{2k}$. Since summands of fully realizable groups are fully realizable and $\D_{2k}$ is not fully realizable, $\D_{4k}$ is not fully realizable.
\end{proof}

\section{Split extensions of certain cyclic groups} \label{sec:cbyc}

In this section, we will determine the fully realizable groups of the form $\C_n \rtimes \C_p$, where $n$ is a positive integer and $p$ is a prime.

\begin{prop}\label{centralcondition}
Let $p$ be prime and let $A$ be an abelian group. Let $\phi \colon \C_p = \langle \sigma \rangle \longrightarrow \text{Aut}(A)$ be a homomorphism such that $\phi_{\sigma}$ belongs to $Z(\text{End}(A))$. If $A \rtimes \C_p$ is fully realizable, then $A$ is fully realizable.
\end{prop}

\begin{proof}
If $\phi$ is trivial, then $A$ is fully realizable by Proposition \ref{prprop} (\ref{PQsd}). So, assume $\phi$ is nontrivial.  Since $Z(\text{End}(A))$ is a subgroup of $\text{End}(A)$, if $\phi_{\sigma}$ belongs to $Z(\text{End}(A))$, then so does $\phi_{\sigma^i}$ for all $i$, showing that  $\text{Im}(\phi) \subseteq Z(\text{End}(A))$. Moreover, $A$ is a normal abelian subgroup of $A \rtimes \C_p$ and the index of $A$ in $A \rtimes \C_p$ is prime. This means $A$ is a maximal abelian subgroup of $A \rtimes \C_p$.  Now the result follows from Proposition \ref{max-ab-semi}.
\end{proof}

\begin{corollary}\label{centralcondition-cor}
Let $p$ be prime. If $\C_n \rtimes \C_p$ is fully realizable, then $n \mid 12$.
\end{corollary}


\begin{thm}
Let $p$ be prime. The group $G = \C_n \rtimes \C_p$ is fully realizable if and only if $G$ is abelian and $p = 2$ with $n \mid 12$ or $p = 3$ with $n \mid 4$.
\end{thm}

\begin{proof} 
The `if' direction follows from Proposition \ref{prprop} (\ref{realizablefab}). So we will prove the `only if' direction.

Suppose  $\C_n \rtimes \C_p$  is fully realizable. By Corollary \ref{centralcondition-cor}, $n \mid 12$. Moreover, $\C_p$ is  a retract of $\C_n \rtimes \C_p$. This means $\C_p$ must be realizable, or equivalently $p =2$ or $3$. The groups to be considered are $\C_2$, $\C_2 \rtimes \C_2$, $\C_3 \rtimes \C_2$, $\C_4 \rtimes \C_2$, $\C_6 \rtimes \C_2$, $\C_{12} \rtimes \C_{2}$, $\C_3$, $\C_2 \rtimes \C_3$, $\C_3 \rtimes \C_3$, $\C_4 \rtimes \C_3$, $\C_6 \rtimes \C_3$, and $\C_{12} \rtimes \C_3$. Among the abelian such semidirect products (which would then be direct products), the fully realizable ones are exactly the groups listed in the statement of the theorem (by Proposition \ref{prprop} (\ref{realizablefab})). So it remains to prove that any such nonabelian semidirect products are not fully realizable.

The automorphism groups $\text{Aut}(\C_3)$, $\text{Aut}(\C_4)$, $\text{Aut}(\C_6)$ and $\text{Aut}(\C_{12})$ have no elements of order 3. Therefore, the corresponding semidirect products $\C_3 \rtimes \C_3$, $\C_4 \rtimes \C_3$ $\C_6 \rtimes \C_3$, and $\C_{12} \rtimes \C_3$ are all abelian. The remaining nonabelian semidirect products are the dihedral groups $\C_3 \rtimes \C_2$, $\C_4 \rtimes \C_2$, $\C_6 \rtimes \C_2$, and $\C_{12} \rtimes \C_{2}$, which we proved are not fully realizable in the previous section.
\end{proof}

\section{Quaternion groups}

For $n \ge 1$, the generalized Quaternion group $\Q_{4n}$ is a group of order $4n$ given by the following presentation:
\[ \Q_{4n} = \langle a, b \,|\, a^{2n} = e, b^4 =1, a^{n}= b^2, bab^{-1}= a^{-1} \rangle.\]
It is an extension of $\C_2$ by $\C_{2n}$:
\[ 1 \xrightarrow[]{} \C_{2n} \xrightarrow[]{} \Q_{4n} \xrightarrow[]{} \C_2 \xrightarrow[]{} 1.\] In this section, we prove the following theorem.

\begin{thm}
The generalized quaternion group $\Q_{4n}$ is fully realizable if and only if $n = 1$.
\end{thm}
\begin{proof} If $n$ is odd then $\mathbf{Q}_{4n} = \C_n \rtimes \C_4$, where the generator of $\C_4$ acts on $\C_n$ by inversion; these groups are not fully realizable for $n > 1$ by Proposition \ref{dicp} below. If $n = 2^k$, then $\Q_{4n}$ is realizable if and only if $n=2$ by \cite{chebolu2019fuchs}, and $\Q_8$ is not fully realizable by Proposition \ref{q8} below. If $n = 2^k m$ where $m > 1$ is odd and $k > 1$, then $\mathbf{Q}_{4n} = \C_m \rtimes \Q_{2^{2+k}}$ is not fully realizable since $\Q_{2^{2+k}}$ is not fully realizable. 
\end{proof}

\begin{prop} \label{dicp}
Let $n$ be odd, write $\C_4 = \langle x \rangle$, and suppose $x$ acts on $\C_n$ via inversion: $xtx^{-1} = t^{-1}$ for all $t \in \C_n$. The semidirect product $G=\C_n \rtimes \C_4$ determined by this action is not realizable in characteristic $2$.
\end{prop}
\n Note that the group appearing in this proposition is the dicyclic group of order $4n$ (this is because $n$ is odd; it's not the dicyclic group when $n$ is even). If $n = 2^k m$ is even with $m$ odd, then $\C_n \rtimes \C_4 = \C_{2^k} \rtimes(\C_m \rtimes \C_4)$, so when $n$ is even and not a power of 2, this proposition implies $\C_n \rtimes \C_4$ is not fully realizable. If $n$ is a power of $2$, then we must consider the group $\C_{2^k} \rtimes \C_4$. We don't have a complete answer in this case, though we do know, for example, that---in contrast to the above proposition---$\C_4 \rtimes \C_4$ is realizable but not fully realizable (\cite[\href{https://cocalc.com/klockrid/main/FPENAG/files/Order-16.ipynb}{Order-16.ipynb}]{ourcode}).
\begin{proof}

%

Assume to the contrary that the ring $R = \F_2[G]/I$ has group of units $G$. Let $y$ denote a generator of $\C_n$. Note that $x^2$ is in the center of $G$ and is the unique element of order 2 in $G$. Since $1 +x^2$ lies in the center of $R$ and its square is zero, the ideal $I = (1 + x^2)$ is contained in the Jacobson radical and $1 + I$ is an elementary abelian normal 2-subgroup of $G$. This forces $I = \{0, 1 + x^2\}$ from which it follows that $g(1 + x^2) = 1 + x^2$ for all $g \in G$.

Consider the element \[\theta = 1 + (y + y^{-1})(1 + x).\] Note that $y + y^{-1}$ is the sum of the elements in the conjugacy class of $y$, so it is in the center of $R$. Hence, \[ \theta^2 =  1 + (y^2 + y^{-2})(1 + x^2) = 1 + (1 + x^2) + (1 + x^2) = 1.\] This forces $\theta = 1$ or $\theta = x^2$. In either case, $\theta$ commutes with $y$, and so $(y + y^{-1})x$ commutes with $y$. By definition of the action of $x$ on $\langle y \rangle$, $xy = y^{-1}x$. Thus:
\[
\begin{aligned}
y(y + y^{-1})x &
= (y + y^{-1})xy \\
(y^2 + 1)x &= (1 + y^{-2})x \\
y^2 + 1 & = y^{-2} + 1 \\ 
y^{2} &= y^{-2} \\
y^4 & = 1.
\end{aligned}
\] This contradicts the fact that $|y|$ is odd.
\end{proof}

\begin{prop} \label{q8}
The quaternion group $\Q_8$ is not fully realizable.
\end{prop}
\begin{proof} Write \[\mathbf{Q}_8 = \langle x, y \, | \, x^4 = 1, x^2 = y^2, yx = x^{-1}y\rangle.\] By Proposition \ref{prprop} (\ref{I-invariant}), if $\Q_8$ is fully realizable, then there is an ideal $I$ in $\F_2[\Q_8]$ that fully realizes $\Q_8$.

In the group algebra $\F_2[\Q_8]$, $\alpha$ is a unit if and only if it has augmentation 1 since $\Q_8$ is a finite 2-group. Hence, $1 + x + y$ is a unit. If $I \subseteq \F_2[\Q_8]$ fully realizes $\Q_8$, then we must have $1 + x + y + u \in I$ for some $u \in \Q_8$. Using GAP in SageMath, we checked that for each possible value of $u$, either the unit group of the quotient ring is too small or there is an endomoprhism of $\Q_8$ that does not preserve $I$ (\cite[\href{https://cocalc.com/klockrid/main/FPENAG/files/Q8.ipynb}{Q8.ipynb}]{ourcode}). For example, if $u = xy$, then $\F_2[\Q_8]/I$ realizes $\Q_8$, but the endomorphism of $\mathbf{Q}_8$ defined by $x \mapsto x$ and $y \mapsto xy$ sends $1 + x + y + xy$ to $1 + x + xy + x^2y$, which is not an element of $I = (1 + x + y + xy)$.
\end{proof}

\section{Finite simple groups}

Recall that a nontrivial group $G$ is said to be simple if its only normal subgroups are the trivial group and the group $G$ itself. In this section, we determine the fully realizable finite simple groups.

The realizable finite simple groups are classified in \cite{davis2014finite}. The authors prove that a finite simple group $G$ is realizable if and only if it is isomorphic to a cyclic group of order $2$,  a cyclic group of prime order $2^k-1$, or a projective special linear group $\mathbf{PSL}_n(\F_2)$ for some $n\ge 3$. Since every fully realizable group is also realizable, we may use their result as a starting point. The next proposition implies that a matrix ring cannot fully realize its group of units.

\begin{prop} \label{transpose-inverse}
Let $n \ge 3$ and let $F$ be a field of characteristic $2$. The automorphism $A \mapsto (A^T)^{-1}$ of  $\GL_n(F)$ is not induced by any ring endomorphism of $\mathbf{M}_n(F)$. 
\end{prop}

\begin{proof} We begin by showing that the identity matrix $I_n$ ($n \ge 3$)  can be written as a sum of two invertible matrices that are not mutually inverse to each other: $ I_n = A + B$ such $AB \ne I_n$. We consider two cases based on the parity of $n$.

If $n$ is odd, write  $n = 2r+3$ ($r \ge 0$), and consider
\[
\begin{bmatrix}
I_r & 0 & 0 \\
0 & I_r & 0\\
0 & 0 & I_3
\end{bmatrix}
= 
\begin{bmatrix}
I_r & I_r & 0 \\
I_r & 0 & 0\\
0  & 0 & U
\end{bmatrix}
+ 
\begin{bmatrix}
0 & I_r & 0 \\
I_r & I_r & 0\\
0 & 0 &  V
\end{bmatrix},
\]
where 
$ U = \begin{bmatrix}
0 & 1 & 1 \\
0 &0 & 1\\
1 & 0 &  0 
\end{bmatrix} $
and 
$ V= \begin{bmatrix}
1 & 1 & 1 \\
0 & 1 & 1\\
1 & 0 &  1
\end{bmatrix}
$. If $n$ is even, write $n = 2r+4$ ($r \ge 0$), and consider
\[
\begin{bmatrix}
I_r & 0 & 0 \\
0 & I_r & 0\\
0 & 0 & I_4
\end{bmatrix}
= 
\begin{bmatrix}
I_r & I_r & 0 \\
I_r & 0 & 0\\
0  & 0 & S
\end{bmatrix}
+ 
\begin{bmatrix}
0 & I_r & 0 \\
I_r & I_r & 0\\
0 & 0 &  T
\end{bmatrix},
\]
where 
$ S = \begin{bmatrix}
1 & 0 & 1 & 1 \\
0 & 1 & 0 & 1\\
1 & 0 & 1 & 0 \\
0 & 1 & 0 & 1
\end{bmatrix} $
and 
$ T = \begin{bmatrix}
0 & 0 & 1 & 1 \\
0 & 0 & 0& 1\\
1 & 0 & 0& 0 \\
0 & 1 & 0 & 0
\end{bmatrix}.
$ It is easy to check that in both cases, we have $A+B = I_n$ but $AB \ne I_n$.

Now, suppose that the given automorphism ($A \mapsto (A^T)^{-1}$ of  $\GL_n(F)$) is induced by a ring endomorphism $\phi$ on $\mathbf{M}_n(F)$ ($n \ge 3$). Applying the homomorphism $\phi$ on both sides of the equation $I_n = A + B$ and letting $C$ and $D$ denote the transposes of $A$ and $B$ respectively gives the following:

\begin{eqnarray*}
I_n  & = & A + B \\
 \phi(I_n) &=& \phi( A + B )\\
 I_n & = & (A^T)^{-1} + (B^T)^{-1} \\
 I_n & = & C^{-1} + D^{-1}\\
 I_n & = & C^{-1} + (I_n - C)^{-1} \\
 I_n - C^{-1} & = & (I_n - C)^{-1} \\
 (I_n - C)(I_n - C^{-1}) & = & I_n\\
 I_n - C - C^{-1} + I_n & =& I_n \\
 I_n - C & = &   C^{-1}\\
 (I_n - C)C & = &  I_n\\
 C^T(I_n - C^T)  & = & I_n\\
 A(I_n -A) & = & I_n\\
 AB & = & I_n.
\end{eqnarray*}
Since $AB \ne I_n$, we have a contradiction, and the proof is complete.
\end{proof}

\begin{thm} \label{main-simple}
A  finite simple group is fully realizable if and only if it is either $\C_2$ or $\C_3$.
\end{thm}
 
\begin{proof}
The only fully realizable cyclic groups of prime order are $\C_2$ and $\C_3$ by Proposition \ref{prprop} (\ref{realizablecyclics}). It therefore suffices to consider the groups $\mathbf{PSL}_n(\F_2) = \GL_n(\F_2)$ for $n \geq 3$.

Suppose $R$ is a ring of characteristic 2, generated by its units, that fully realizes its unit group $G = \GL_n(\F_2)$ for some $n \geq 3$. Since $R$ is finite, the Jacobson radical $J$ of $R$ is nilpotent. If $J$ is nonzero, then there is a nonzero ideal $I$ of $R$ whose square is zero. Hence, $1 + I$ is a nontrivial elementary abelian normal subgroup of $G$. Since $G$ has no such subgroup, we conclude that $J = 0$ and $R$ is a product of matrix rings over fields of characteristic 2 by the Artin-Wedderburn theorem. Since $R^\times = \GL_n(\F_2)$ is simple and indecomposable, $R$ must be either $\mathbf{M}_n(\F_2)$ or $\mathbf{M}_n(\F_2) \times \F_2$. Note that we cannot have more than one factor of $\F_2$ because if we did then the ring would not be generated by its units.
 
 If $R=\mathbf{M}_n(\F_2)$, Proposition \ref{transpose-inverse} shows that $R$ does not fully realize $\GL_n(\F_2)$.
 If $R = \mathbf{M}_n(\F_2) \times \F_2$, then the endomorphisms of $R$ are the maps of the form $(A, x) \mapsto (\phi(A), x)$ and $(A, x) \mapsto  (xI_n, x)$ where $\phi$ is an endomorphism of $\mathbf{M}_n(\F_2)$. (This is because $\mathbf{M}_n(\F_2)$ and $\F_2$ do not have any nontrivial ideals.) Therefore, we will be done if we can produce an endomorphism of $\GL_n(\F_2)$ that is not realized by the ring $\mathbf{M}_n(\F_2)$ when $n \ge 3$. This is achieved by invoking Proposition \ref{transpose-inverse} again. This completes the proof of the theorem.
\end{proof}

\section{Symmetric and alternating groups}
For any positive integer $n$, let $\mathbf{S}_n$ denote the symmetric group on $n$ letters.  In \cite{davis2014alternating} the authors prove that $\mathbf{S}_n$ is realizable if and only if  $n \le 4$. Note that $\mathbf{S}_1 = \C_1$ and $\mathbf{S}_2 = \C_2$ are fully realizable.  The group $\mathbf{S}_3$ is isomorphic to $\D_6$, which is not fully realizable by Theorem \ref{dihedral}. Since $\mathbf{S}_4 = \C_2^2 \rtimes \mathbf{S}_3$ and $\mathbf{S}_3$ is not fully realizable, $\mathbf{S}_4$ is not fully realizable by Proposition \ref{prprop} (\ref{PQsd}). These observations give the following theorem.

\begin{thm} The following are equivalent.
\begin{enumerate}
    \item $n = 1$ or $2$
    \item $\mathbf{S}_n$ is fully realizable 
    \item $\mathbf{S}_n$ is simple.
\end{enumerate}
\end{thm}

Recall that the holomorph $\mathrm{Hol}(G)$ of a group $G$ is constructed using the natural action of the automorphism group of $G$ on $G$: $\mathrm{Hol}(G) \cong G \rtimes \text{Aut}(G)$. Since $\text{Aut}(G)$ is a retract of $\text{Hol}(G)$, and for $n \ne 2$ or $6$ $\text{Aut}(\mathbf{S}_n) \cong \mathbf{S}_n$, we obtain the following corollary.

\begin{corollary}
If $n \ne 2$ or $6$, $\text{Hol}(\mathbf{S}_n)$ is not fully realizable.
\end{corollary}

Next, we consider the alternating groups $\A_n$ for $n >1$. Recall that $\A_n$ is the subgroup of even permutations in $\mathbf{S}_n$ and has order $n!/2$. The groups $\A_2 = \C_1$ and $\A_3 = \C_3$ are fully realizable. For $n \geq 5$, $\A_n$ is simple and therefore not fully realizable by Theorem \ref{main-simple}. The group $\A_4 = \C_2^2 \rtimes \C_3$ is fully realizable by \cite[\href{https://cocalc.com/klockrid/main/FPENAG/files/A4.ipynb}{A4.ipynb}]{ourcode}. 

\begin{thm}
The following are equivalent.
\begin{enumerate}
    \item $n=4$
    \item $\A_n$ is fully realizable
    \item $\A_n$ is not simple
\end{enumerate}
\end{thm}

\section{Finite general linear groups} \label{sec:fglg}

In this section we will determine the general linear groups over finite fields ($\GL_n(\F_q$)) that are fully realizable.  We first recall two facts we will use in our analysis.

\begin{prop}
For any finite ring $R$, let $J$ denote its Jacobson radical. Then $1+J$ is a normal subgroup of $R^\times$ such that \[ (R/J)^\times \cong R^\times/(1+J).\]
Moreover, $|1+J|$ divides $|R|$. In particular, if $R$ has characteristic $2$, then $1+J$ is a normal $2$-subgroup of $R^\times$.
\end{prop}

\begin{thm}[Dickson]  Every normal subgroup of $\GL_n(\F_q)$ either contains $\SL_n(\F_q)$ or is contained in the center $\C_{q-1}$ with the exception of $n=2$ and $q=2, 3$.
\end{thm}

\begin{thm} The group $\GL_n(\F_q)$ is fully realizable if and only if $n=1$ and $q$ is $2 ,3, 4, 5, 7$ or $13$.
\end{thm}

\begin{proof}
If $n=1$, then $\GL_1(\F_q) = \C_{q-1}$. This group is fully realizable if and only if $q-1$ divides 12. This proves one direction of our theorem.

For the other direction, assume $n\ge 2$, and let $q$ be any prime power.  By Example \ref{generallinear}, if $\GL_n(\F_q)$ is fully realizable, then $q-1$ must divide $12$, which means $q$ is $2 ,3, 4, 5, 7$ or $13$.

The group $\GL_n(\F_2)$ is isomorphic to $\mathbf{S}_3$ when $n=2$, and it is a finite simple group for $n > 2$. We have shown that none of these groups are fully realizable.

The group $\GL_2(\F_4)$ has $\A_5$ as a summand, so it cannot be fully realizable. For $n \ge 3$, we will show that  $\GL_n(\F_4)$ cannot be  fully realizable.
Suppose $R$ is a ring of characteristic 2, generated by its units, that fully realizes its unit group $G = \GL_n(\F_4)$ for some $n \geq 3$. Since $R$ is finite, the Jacobson radical $J$ of $R$ is nilpotent. If $J$ is nonzero, then there is a nonzero ideal $I$ of $R$ whose square is zero. Hence, $1 + I$ is a nontrivial elementary abelian normal subgroup of $G$. Since $G$ has no such subgroup, we conclude that $J = 0$ and $R$ is a product of matrix rings over fields of characteristic 2 by the Artin-Wedderburn theorem. Note that $R^\times = \GL_n(\F_4)$ has only two non-trivial normal subgroups ($\C_3$ and $\SL_n(\F_4)$), has center $\C_3$, and $R$ is assumed to be generated by its units. These facts force $R$ to be either $\mathbf{M}_n(\F_4)$ or $\mathbf{M}_n(\F_4) \times \F_2$. Using the argument in the last paragraph of the proof of Proposition \ref{main-simple}, neither of these rings fully realizes its group of units.

We will next show that $\GL_n(\F_q)$  is not fully realizable when $q$ is $3, 5, 7,$ or $13$ and $n \ge 2$. Since any fully realizable group is realizable in characteristic $2$, we  will be done if we can  show that these groups are not realizable in characteristic $2$. 
Assume to the contrary that $\GL_n(\F_q)$ is realizable by a ring $R$ in characteristic $2$. Without loss of generality, we may assume that $R$ is a finite ring that is generated by its units.  Let $J$ denote the Jacobson radical of $R$.  Then $R/J$ is a semisimple ring whose unit group is isomorphic to  $R^\times/(1+J)$.  Since $1+J$ is a normal $2$-subgroup of $\GL_n(\F_q)$, it must be contained in the center because $\SL_n(\F_q)$ is not a $2$-group. This shows that $1+J \cong \C_{2^k}$ for $0 \le k \le 2$; these are the only cyclic 2-groups that occur as subgroups in one of the multiplicative groups  $(\F_q)^\times$. Using the Artin-Wedderburn theorem, we have:
\[\GL_n(\F_q)/\C_{2^k} \cong \prod_i \GL_{n_i}(\F_{2^{r_i}}) \] 
The order of the center of the right hand side is odd, so $\C_{2^k}$ is the Sylow 2-subgroup of $\C_{q-1} = Z(\GL_n(\F_q))$. 

If $q = 3$ or $5$, then the center of the left hand side is trivial and thus $r_i = 1$ for all $i$. If $q = 7$ or $13$, then the center of the left hand side is $\C_3$ and thus $r_i = 1$ for all $i$ except for one value of $i$ where we have $r_i = 2$. We will next consider two cases, depending on whether $\GL_n(\F_q)/\C_{2^k}$ is indecomposable. Before doing so, we need a lemma.

\begin{lem} Let $v \geq 2$ be an integer and let $t$ be a prime power where $t > 3$ if $v = 2$. Then, $\PGL_v(\F_t)$ is simple if and only if $\gcd(v, t-1) = 1$ if and only if $\PGL_v(\F_t) = \PSL_v(\F_t)$.
\end{lem}
\begin{proof} Let $Z = \F_t^\times = \C_{t-1}$ denote the center of $\GL_v(\F_t)$. If $\PGL_v(\F_t)$ is simple, then we must have that the normal subgroup $\langle \SL_v(\F_t), Z\rangle$ of $\GL_v(\F_t)$ is equal to $\GL_v(\F_t)$. Hence the determinant map must induce an isomorphism from $Z$ to $\C_{t-1}$; this is the $v$-th power map, so $\gcd(v, t-1) = 1$. If $\gcd(v, t-1) = 1$, then $\PSL_v(\F_t) = \SL_v(\F_t)$, and the natural map (inclusion followed by the quotient map) \[\SL_v(\F_t) \longrightarrow \GL_v(\F_t) \longrightarrow \PGL_v(\F_t)\] is injective and therefore an isomorphism since the source and the target have the same order. If  $\PGL_v(\F_t) = \PSL_v(\F_t)$, then $\PGL_v(\F_t)$ is simple since $\PSL_v(\F_t)$ is simple.
\end{proof}

Suppose $\GL_n(\F_q)/\C_{2^k}$ is indecomposable. Then, it is isomorphic to $\GL_m(\F_2)$ for some $m \geq 2$ or $\GL_m(\F_4)$ for some $m \geq 2$. The order of $\GL_n(\F_q)$ is at least 48 and the order of $\GL_2(\F_2)$ is 6, so we may assume $m > 2$ in the first of these two cases. 

Suppose $\GL_n(\F_q)/\C_{2^k} = \GL_m(\F_2)$, in which case $q = 3$ or $5$ and we in fact obtain $\PGL_n(\F_q) = \GL_m(\F_2)$. Since $\GL_m(\F_2) = \PSL_m(\F_2)$ is simple ($m > 2$), $\PGL_n(\F_q) = \PSL_n(\F_2)$ and thus $\PSL_n(\F_q) = \PSL_m(\F_2)$. By \cite[\S 1.2]{wilson}, this forces $m = 3, n = 2, q = 7$, but $q = 3$ or $5$, a contradiction.

Suppose $\GL_n(\F_q)/\C_{2^k} = \GL_m(\F_4)$, in which case $q = 7$ or $13$. Modding out by centers, we obtain $\PGL_n(\F_q) = \PGL_m(\F_4)$. If $3 \nmid m$, then $\PGL_m(\F_4) = \PSL_m(\F_4)$ (a simple group) and so $\PSL_n(\F_q) = \PSL_m(\F_4)$; by \cite[\S 1.2]{wilson}, $n = m = 2$ and $q = 5$, but $q$ is 7 or 13, a contradiction. If $3 \mid m$, then the normal subgroups of $\GL_m(\F_4)$ are \[ \{1 \} \leq \C_3 \leq \SL_m(\F_4) \leq \GL_m (\F_4).\] Since $\GL_n(\F_q)/\C_{2^k} = \GL_m(\F_4)$, the normal subgroups of $\GL_n(\F_q)/\C_{2^k}$ are \[\{1\} \leq \C_3 = \C_{q-1}/\C_{2^k} \leq \langle \SL_n(\F_q), \C_{2^k}\rangle/\C_{2^k} \leq \GL_n(\F_q)/\C_{2^k}.\] From this it follows that $\SL_m(\F_4) = \SL_n(\F_q)$. Modding out by centers, we obtain $\PSL_m(\F_4) = \PSL_n(\F_q)$, which again by \cite{wilson} implies $q = 5$, a contradiction.

Next we turn to the case where $\GL_n(\F_q)/\C_{2^k}$ is decomposable. In this case, there must be two nontrivial normal subgroups $N, N'$ of $\GL_n(\F_q)$ whose intersection is $\C_{2^k}$ and whose product is the whole group. This means both cannot contain $\SL_n(\F_q)$, and thus one must be precisely equal to $Z(\GL_n(\F_q)) = \C_{q-1}$ since the index of $\C_{2^k}$ in  $\C_{q-1}$ must be 3 (if $q = 3$ or $5$, then we cannot have two distinct nontrivial normal subgroups whose intersection is $\C_{2^k} = Z$). Thus, one of the direct factors is $\C_3$ and we have \[\GL_n(\F_q)/\C_{2^k} = \GL_m(\F_2) \times \C_3.\] Modding out by centers on both sides, we obtain \[\PGL_n(\F_q) = \PGL_m(\F_2).\] This forces (recall that $m > 2$, and the right hand side is simple) \[\PSL_n(\F_q) = \PSL_m(\F_2).\] By \cite{wilson}, this means $n = 2, q = 7, m = 3$ and $k = 1$. But $\PGL_2(\F_7)$ is not equal to $\PSL_2(\F_7)$ (since $\gcd(2, 7-1) = 2 \neq 1$), a contradiction. This completes the proof.
\end{proof}

\section{Groups with a cyclic Sylow 2-subgroup}

It was shown in \cite{chebolu2022fuchs} that a finite group of odd order is fully realizable if and only if it is either $\C_1$ or $\C_3.$ We therefore restrict ourselves to groups of even order. Such a group has a nontrivial Sylow $2$-subgroup. We now classify the fully realizable groups of even order with a cyclic Sylow 2-subgroup.

\begin{prop} \label{Sylow2iscyclic}
Let $G$ be a finite group of even order that is fully realizable, and let $H$ be the Sylow $2$-subgroup of $G$. If $H$ is cyclic,  then  $G$ is either $L \rtimes \C_2$ or $L \rtimes \C_4$ where $L$ is a group of odd order. 
\end{prop}

\begin{proof}
Since the Sylow $2$-subgroup $H$ is cyclic, by a theorem of Cayley, we know that $H$ has a normal complement $L$. Now $L$ is subgroup of $G$ with odd order and $G = HL$ with $H \cap L = \{ e \}$. This shows that  $H$ is a retract of $G$. Since $G$ is fully realizable, so is $H$. The only powers of $2$ that divide $12$ are $2$ and $4$.  This means $H$ must be $\C_2$ or $\C_4$, and $G$ is either $L \rtimes \C_2$ or $L \rtimes \C_4$. 
\end{proof}

We now wish to determine the possible values of $L$. It is useful to first consider general linear groups of order $2s$ and $4s$ where $s$ is odd.


\begin{prop}\label{glp}
Let $s$ be an odd positive integer and let $G = \GL_m(\F_{2^k})$. If $G$ has order $2s$, then $G=\mathbf{S}_3 $. If  $G$ has order $4s$, then $G = \C_3 \times \mathbf{A}_5$.
\end{prop}
\begin{proof}
The exponent of 2 in \[|\GL_m(\F_{2^k})| = (2^{mk} -1)(2^{mk}-2^k) \cdots (2^{mk}-2^{(m-1)k})\] is $km(m-1)/2$. Hence, we must have $k = 1$ and $m = 2$ or $k = 2$ and $m = 2.$ Now $\GL_2(\F_2) = \mathbf{S}_3$ and $\GL_2(\F_4) = \C_3 \times \mathbf{A}_5.$
\end{proof}

\begin{prop} \label{lc2}
Let $L$ be a group of odd order. The group $G = L \rtimes \C_2$ is fully realizable if and only if  it  is isomorphic to either $\C_2$ or  $\C_3 \times \C_2$.
\end{prop}

\begin{proof} The `if' direction is clear. For the converse, suppose that  $G$ is fully realized by a finite ring $R$ of characteristic 2. We have two cases to consider.

Suppose the Jacobson radical of $R$ is trivial. Then, $R$ is semisimple, and by the Artin-Wedderburn theorem we can write $R \cong \prod \mathbf{M}_{n_i}(\F_{2^{r_i}})$. Since the unit group of every factor must be fully realizable and $|R^\times| = 2|L|$, Proposition \ref{glp} implies that $n_i = 1$ for all $i$ ($\mathbf{S}_3$ is not fully realizable) and hence $G$ is abelian, forcing $L = \C_1$ or $L = \C_3$.

Suppose the Jacobson radical $J$ of $R$ is nontrivial. Since $R$ is artinian, $J$ is nilpotent, and there is a nontrivial ideal $I$ (an appropriate power of $J$) such that $I^2 = 0$. Then, $1+I$ is normal in $L \rtimes \C_2$  (see \cite[2.6]{chebolu2017fuchs}) and is an elementary abelian 2-group since $(1+x)^2 = 1$ for all $x \in I$. This implies that $1+I = \C_2$. Now both $L$ and $\C_2$ are normal in $G$ and thus $G = L \times \C_2$. Since direct factors of fully realizable groups are fully realizable, $L$ must be fully realizable. Since $L$ is a group of odd order, it must be either $\C_1$ or $\C_3$.
\end{proof}

\begin{prop} \label{Lab}
Let $L$ be a finite group of odd order. If $L\rtimes \C_4$ is fully realizable, then $L$ is abelian.
\end{prop}
\begin{proof}

Suppose $R$ is a ring a characteristic 2 that fully realizes $G = L \rtimes \C_4$ and let $S$ be the subring of $R$ generated by $L$. Since $S$ is a quotient of the semisimple ring $\F_2[L]$, it is a product of matrix rings of the form $\GL_m(\F_{2^k})$ by the Artin-Wedderburn theorem. However, we cannot have $S^\times = G$ since neither $\mathbf{A}_5$ nor $\mathbf{S}_3$ is fully realizable (see Proposition \ref{glp}). So $S = F \times \GL_2(\F_2)^\epsilon$ where $F$ is a finite product of fields and $\epsilon \in \{0, 1\}$. Since $L$ is a subgroup of $S^\times = F^\times \times \mathbf{S}_3$ of odd order, it must be abelian, since $F^\times \times \mathbf{S}_3$ has no nonabelian subgroup of odd order.
\end{proof}

\begin{prop} Let $L$ be a group of odd order. The group $G = L \rtimes \C_4$ is fully realizable if and only if  it  is isomorphic to either $\C_4$ or  $\C_3 \times \C_4$.
\end{prop}
\begin{proof}
The `if' direction is clear. For the converse, let $R$ be a ring of characteristic 2 that fully realizes $G$. Write $\C_4 = \langle x \rangle$. If $R$ is semisimple, then as in the proof of Proposition \ref{lc2}, we obtain that $G$ is abelian, forcing $G = \C_4$ or $G = \C_3 \times \C_4$.

If $R$ is not semisimple, then we may again follow the argument in Proposition \ref{lc2} to conclude that $\langle x^2 \rangle$ is normal in $G$; this implies that $x$ acts on $L$ via an automorphism of order 2. By Proposition \ref{Lab}, $L$ is abelian. By Proposition \ref{fpfsplit}, $G = H_f \times (H' \rtimes\C_4)$ where the action of $x$ on $H'$ is fixed-point free. Examining the proof of Proposition \ref{fpfsplit}, we see that $H' = \oplus_q B_q$, where the indices $q \neq 2$ are primes and $B_q$ is a module over $\Z_{q^k}$ on which $x$ acts via inversion. The action of $x$ on $B_q$ therefore acts decomposably on any direct sum decomposition of $B_q$, and choosing any such summand $\C_{q^l}$, $G$ has the form $G = H_f \times (H'' \rtimes (\C_{p^l} \rtimes \C_4))$. By Proposition \ref{prprop} (\ref{PQsd}), $\C_{p^l} \rtimes \C_4$ is fully realizable. This contradicts Proposition \ref{dicp} unless $H' = 0$. Since $H' = 0$, $G$ is abelian, and the proof is complete.
\end{proof}

By combining the above propositions we get the follow theorem.

\begin{thm}
The fully realizable finite groups of even order with a cyclic Sylow $2$-subgroup are $\C_2$, $\C_6$, $\C_4$ and $\C_{12}$.
\end{thm}

\section{Split extensions of elementary abelian 2-groups by cyclic $p$-groups}

In this section, we consider certain groups of the form $G = A \rtimes \C_{p^n},$ where $A$ is a finite 2-group and $p^n$ is a prime power. We construct three infinite families of fully realizable nonabelian groups. Note that if $G$ is fully realizable, then $\C_{p^n}$ is fully realizable by Proposition \ref{prprop} (\ref{PQsd}). Thus, by Proposition \ref{prprop} (\ref{realizablecyclics}), $p^n = 1, 2, 3,$ or $4$. Our first goal is to prove the following theorem, characterizing the fully realizable groups of the form $A \rtimes \C_{p^n}$ when $A$ is a finite elementary abelian 2-group. In a semidirect product $X \rtimes G$ determined by a homomorphism $\phi\colon G \longrightarrow \mathrm{Aut}(X),$ if $\phi$ is injective then we say $G$ acts {\bf faithfully} on $X.$

\begin{thm} Let $p^n$ be a prime power, let $A$ be a finite elementary abelian $2$-group, and let $G = A \rtimes \C_{p^n}$. If $G$ is fully realizable, then $p^n = 1, 2, 4$ or $3$. \label{ea2gpn}
\begin{enumerate}
    \item If $p^n = 2$, then $G$ is fully realizable if and only if $G$ is abelian.
    \item If $p^n = 4$, then $G$ is fully realizable if and only if $\C_4$ does not act faithfully on $A$.
    \item If $p^n = 1$ or $3$, then $G$ is fully realizable.
\end{enumerate}
\end{thm}

\n Note that when $A$ is a finite elementary abelian 2-group and $p^n = 1$, we have $G = A,$ and $A$ is fully realizable by \cite[5.1]{chebolu2022fuchs}.

Our approach is as follows. In \S \ref{sec:general}, we prove a number of preliminary results to help establish that certain specific groups are fully realizable. In \S \ref{sec:n2}, \S \ref{sec:n4} and \S \ref{sec:n3}, we consider $p^n = 2, 4,$ and $3,$ respectively. In each case, we reduce the question of whether $A \rtimes \C_{p^n}$ is fully realizable to whether $M \rtimes \C_{p^n}$ is fully realizable for specific indecomposable $\F_2[\C_{p^n}]$-modules $M$. We then use GAP to settle these specific cases.

The case $p = 3$ is special and may be generalized; see Theorem \ref{na2sg} in \S \ref{sec:n3}. We are also able to construct a family of fully realizable groups of the form $A \rtimes \C_3$ where $A$ is not abelian in \S \ref{sec:nab2}. The groups in this family have the form $A_0 \rtimes \C_6 = (A_0 \rtimes \C_2) \rtimes \C_3,$ where $A_0$ is a finite elementary abelian $2$-group. As above, we study indecomposable $\F_2[\C_6]$-modules $M$ and use GAP to determine whether the groups $M\rtimes \C_6$ are fully realizable.

\subsection{Tools} \label{sec:general} The following proposition enables us, in many cases, to identify the correct $\mathrm{End}(G)$-invariant ideal $I$ that fully realizes $G$.

\begin{prop} \label{centralizercriterion}
Let $V$ and $W$ be finite groups with $V$ abelian, let $G = V \rtimes W$, and let $I$ denote the ideal in $\F_2[G]$ generated by all elements of the form \[1 + s + t + st\] where $s \in C_G(V)$ and $t \in G$. 
\begin{enumerate}
    \item If every endomorphism $\rho$ of $G$ satisfies $\rho(C_G(V)) \subseteq C_G(V)$, then $I$ is an $\mathrm{End}(G)$-invariant ideal. 
    \item If $G$ embeds in $\F_2[G]/I$ and $I$ is $\mathrm{End}(G)$-invariant, then  $\rho(C_G(V)) \subseteq C_G(V)$ for all $\rho \in \mathrm{End}(G)$.
\end{enumerate}
\end{prop}
\begin{proof}
The first claim is true by definition of $I$. We prove the contrapositive of the second claim. Suppose there is some $\rho \in \mathrm{End}(G)$ and $s \in C_G(V)$ such that $\rho(s) \not \in C_G(V)$ and $I$ is $\mathrm{End}(G)$-invariant. Then, there is an element $v_0 \in V$ such that $\rho(s)v_0 \neq v_0\rho(s)$. Since $1 + s + \rho^{-1}(v_0) + s\rho^{-1}(v_0) \in I,$ we have $1 + \rho(s) + v_0 + \rho(s)v_0 \in I$ by the $\mathrm{End}(G)$-invariance of $I$. On the other hand, since $V$ is abelian, $v_0 \in C_G(V)$, and $1 + v_0 + \rho(s) + v_0\rho(s) \in I$. Hence, the distinct elements of $G$, $\rho(s)v_0$ and $v_0\rho(s)$, are congruent modulo $I$, and thus $G$ does not embed in $\F_2[G]/I$.
\end{proof}

It is obvious that if $G$ and $H$ are groups realizable in characteristic 2, then so is their product. However, the next proposition identifies a specific ideal of $\F_2[G\times H]$ that realizes the product.

\begin{prop} \label{realizeproduct}
Suppose $\F_2[G]/J$ realizes $G$ and $\F_2[H]/K$ realizes $H$. Then, $\F_2[G \times H]/I$ realizes $G \times H$, where $I$ is generated by $J, K$, and $T = \{1 + g + h + gh \, | \, g \in G, h \in H\}.$
\end{prop}
\begin{proof}
By \cite[3.8]{chebolu2022fuchs}, $S = \F_2[G \times H]/(T)$ is isomorphic to a subring of $\F_2[G] \times \F_2[H]$ that realizes $\F_2[G]^\times \times \F_2[H]^\times$. We then have $S/((J \times K) \cap S) = \F_2[G \times H]/I$ has unit group $G \times H$.
\end{proof}

Let $H = \C_2^L \times \C_4^\epsilon = \oplus_{\alpha \in L} \C_2 \times \C_4^\epsilon,$ where $\epsilon = 0$ or $1$. In \cite[5.6]{chebolu2022fuchs}, we proved that $\F_2[H]/K$ fully realizes $H$, where $K$ is the ideal generated by the elements $1 + u + v + uv$ where $u, v \in H$ and at most one of $u, v$ has order 4. The next proposition says that, for certain fully realizable groups $G$, $G \times H$ is also fully realizable.

\begin{prop}
Let $H = \C_2^L \times \C_4^\epsilon$ where $\epsilon \in \{0,1\}$ and let $K$ be the ideal of $\F_2[H]$ generated by the elements $1 + u + v + uv$ where $u, v \in H$ and at most one of $u, v$ has order 4. Suppose the ideal $J \subseteq \F_2[G]$ fully realizes $G$ where \label{sdpxh}
\begin{enumerate}
\item \label{a0} when $\epsilon = 1$, $G$ has no element of order $4$ and $G$ does not have $\C_4$ as a quotient,
    \item $J$ is generated by elements of the form $1 + x + y + xy$ (where $x, y \in G$), and \label{a1}
    \item for all $x, y \in G$ such that $xy = yx$ and $|x| = 2$, $1 + x + y + xy \in J.$\label{a2}
\end{enumerate}
Define $I$ as in Proposition \ref{realizeproduct}: \[ I = (J, K, 1 + g + h + gh \, | \, g \in G, h \in H) \subseteq \F_2[G\times H].\] The ring $\F_2[G \times H]/I$ fully realizes $G \times H$.
\end{prop}
\n In assumption (\ref{a2}), one can show that under the given conditions on $x$ and $y$ there is unit $u \in \C_G(\langle x, y \rangle)$ with $|u| \mid 2$, $u \neq x$, and  $u + x + y + xy \in J$. So here we suppose $u = 1$.
\begin{proof}
We already know that the indicated ring realizes $G \times H$ by the previous proposition, so we need only prove that the ideal $I$ is $\mathrm{End}(G \times H)$-invariant.

The endomorphisms $F$ of $G \times H$ are in one-to-one correspondence with the tuples $(\phi, \rho, \psi, \tau)$, where $\phi \in \mathrm{End}(G), \tau \in \mathrm{End}(H),$ $\rho \in \mathrm{Hom}(H, G)$, $\psi \in \mathrm{Hom}(G, H)$, and $\rho(x) \in \C_G(\mathrm{Im}\, \phi)$ for all $x \in H$. Note that $F(g) = \phi(g)\psi(g)$ for $g \in G$ and $F(h) = \rho(h) \tau(h)$ for $h \in H$.

We must show that $F$ maps every generator of $I$ into $I$. First, take a generator $\theta$ of $J$, which has the form $\theta = 1 + x + y + xy$ for some $x, y \in G$ by (\ref{a1}). Now, modulo $I$, we have
\[
\begin{aligned}
F(\theta) &= 1 + \phi(x)\psi(x) + \phi(y)\psi(y) + \phi(x)\phi(y)\psi(x)\psi(y) \\
&= (1 + \phi(x) + \phi(y) + \phi(x)\phi(y)) + (1 + \psi(x) + \psi(y) + \psi(x)\psi(y)).
\end{aligned}
\]
The first quantity above lies in $J$, and hence in $I$, since $J$ is $\mathrm{End}(G)$-invariant. The second quantity lies in $K$ since $G$ does not have $\C_4$ as a quotient (assumption (\ref{a0})), and thus $\psi(x)$ and $\psi(y)$ have order at most 2. So $F(\theta) = 0$ modulo $I$.

Next, consider $\theta = 1 + u_1 + u_2 + u_1u_2 \in K$. We have (modulo $I$):
\[
\begin{aligned}
F(\theta) &= 1 + \rho(u_1)\tau(u_1) + \rho(u_2)\tau(u_2) + \rho(u_1)\rho(u_2)\tau(u_1)\tau(u_2) \\
&= (1 + \rho(u_1) + \rho(u_2) + \rho(u_1)\rho(u_2)) + (1 + \tau(u_1) + \tau(u_2) + \tau(u_1)\tau(u_2)).
\end{aligned}
\] The first term in parentheses lies in $J$ by assumption (\ref{a2}) since $\rho(u_1), \rho(u_2)$ are commuting elements, one of which has order at most 2 because one of $u_1, u_2$ has order at most 2. The second term lies in $K$ (at most of of $u_1, u_2$ has order 4, so the same is true for $\tau(u_1), \tau(u_2)$). Hence, $F(\theta) = 0$ modulo $I$.

Finally, consider $\theta = 1 + g + h + gh$, where $g \in G$ and $h \in H$. Modulo $I$, 
\[
\begin{aligned}
F(\theta) &= 1 + \phi(g)\psi(g) + \rho(h)\tau(h) + \phi(g)\rho(h)\psi(g)\tau(h) \\
&= (1 + \phi(g) + \rho(h) + \phi(g)\rho(h)) + (1 + \psi(g) + \tau(h) + \psi(g)\tau(h)).
\end{aligned}
\] The first term lies in $J$ by assumption (\ref{a2}) since $|\rho(h)| \leq 2$ ($G$ has no element of order 4 by condition (\ref{a0})) and $\rho(h)$ commutes with $\phi(g)$ (see the second paragraph of this proof). The second term lies in $K$ since no element in the image of $\psi$ can have order 4. Thus, $F(\theta) = 0$ modulo $I$. This completes the proof.
\end{proof}

Before turning to the case $p^n = 2$, we prove one last proposition concerning the centralizer of the defining normal subgroup in a semidirect product.

\begin{prop}
Let $G = A \rtimes H$, where $A$ is abelian and the semidirect product is determined by the homomorphism $\phi \colon H \longrightarrow \mathrm{Aut}(A).$ Then, $C_G(A) = A \times \ker \phi$. \label{centralizerofSDP}
\end{prop}
\begin{proof}
Take $a \in A$ and $h \in H$. We have
\[
\begin{aligned}
ah \in C_G(A) \iff& aha_0 = a_0 ah \text{ for all $a_0 \in A$} \\
\iff& ha_0 = a_0 h \text{ for all $a_0 \in A$ (since $A$ is abelian)}\\
\iff& ha_oh^{-1} = a_0 \text{ for all $a_0 \in A$} \\
\iff& \phi(h) = 1.
\end{aligned}
\] Thus, $C_G(A) = A \times \ker \phi$.
\end{proof}

\subsection[Case: 2]{Cyclic of order 2} \label{sec:n2}

In this case, the finite elementary abelian 2-group $A$ is a module over $Q = \F_2[\C_2] = \F_2[z]/(z^2)$ (the element $1 + z$ generates $\C_2$). Every finitely generated module over $Q$ is a direct sum of copies of the module $\F_2$, on which $1 + z$ acts trivially, and $Q$ itself. Under the action of $1 + z$ on $Q$: \[0 \mapsto 0, z \mapsto z, 1 + z \mapsto 1, 1 \mapsto 1 + z.\] Switching to multiplicative notation, let $x = 1 + z$, $a = 1$, and $b = z$. The action of $\C_2 = \langle x \rangle$ on $Q$ is summarized by
\[ 
\begin{aligned}
Q &= \C_2 \times \C_2 = \langle a \rangle \times \langle b \rangle &&
a \mapsto a^x = ab,
b \mapsto b^x = b.
\end{aligned}
\]

It follows that we may write $A = A' \times A''$, where $\C_2$ acts decomposably on $A' \times A''$, trivially on $A'$, and $A'' = \oplus_{i=1}^m Q$. Thus, \[A \rtimes \C_2 = A' \times (A'' \rtimes \C_2),\] where $A'$ is a finite abelian 2-group.

\begin{thm}
Let $A$ be a finite abelian $2$-group equipped with a $\C_2$-action. The group $A \rtimes \C_2$ is fully realizable if and only if $\C_2$ acts trivially.
\end{thm}
\n Though $A\rtimes \C_2$ is not fully realizable in this case, it is realizable; we will not give a proof here, as it is nearly identical to the proof given in Section \ref{sec:n4} below that $A\rtimes \C_4$, with $\C_4$ acting non-faithfully, is realizable. The group with GAP id (32, 27), which has the form $\C_2^4 \rtimes \C_2$, is in this category.
\begin{proof}
By the discussion before the statement of the proposition, it suffices to prove that $G = A'' \rtimes \C_2 = \oplus_{i=1}^m Q \rtimes \C_2$ is not fully realizable when $m > 0$. Since $\C_2$ acts decomposably on $\oplus_{i=1}^m Q,$ if $G$ is fully realizable then so is $Q \rtimes \C_2$. This group is the dihedral group of order 8, which we proved is not fully realizable in Section \ref{sec:dihedral}.
\end{proof}

\subsection[Case: 4]{Cyclic of order 4} \label{sec:n4}

If $A$ is a finite elementary abelian $2$-group equipped with a $\C_4$-action, then $A$ is a module over \[S = \F_2[\C_4] = \F_2[z]/(z^4)\] where $1+z$ generates $\C_4$. Every finitely generated $S$-module is a direct sum of copies of $S, \F_2 = S/(z), Q = S/(z^2)$, and $U = S/(z^3)$. The group $\C_4 = \langle x \rangle$ acts trivially on $\F_2$, it acts via $\C_2$ on $Q$ (i.e., $x^2$ acts trivially), and it acts faithfully on $S$ and $S/(z^3)$.

First, we consider the case where the action of $\C_4$ on $A$ is not faithful, so that $A$ is a sum of copies of $\F_2$ and $Q$. Write $A = A' \times A''$  where $\C_4$ acts trivially on $A'$, $A'' = \oplus_{i= 1}^m Q$, and $A \rtimes \C_4 = A' \times (A'' \rtimes \C_4)$. Let $\{a_i, b_i\}$ be a basis for the $i$th copy of $Q$ in $A''$ where $a_i^x = a_i b_i$ and $b_i^x = b_i$. Note that in the semidirect product, this action is conjugation: $t^x = xtx^{-1}$.

\begin{prop} \label{hc4}
Let $I$ denote the two-sided ideal in $\F_2[A'' \rtimes \C_4]$ generated by the elements $1 + a_i + x + a_ix$, $1 + a_i + a_j + a_ia_j$, and $1 + x + x^2 + x^3$ for $i, j \in \{1, \dots, n\}$. The ring $\F_2[A'' \rtimes \C_4]/I$ fully realizes $A'' \rtimes G$.
\end{prop}
\n In the proof, we will show that $I$ is in fact the ideal defined in Proposition \ref{centralizercriterion}. However, we defined it here using what we expect to be a minimal set of generators; this makes it easier to work with such groups in GAP.
\begin{proof}
Let $G = A'' \rtimes \C_4$, let $R = \F_2[G]$, and define $I$ as in the statement of the theorem. One may check that, for all $i$ and $j$, the following expressions lie in $I$, in addition to the generators given in the statement of the theorem: \[
\begin{aligned}
1 + a_i + b_j + a_ib_j &= a_jb_j(1 + a_i + x + a_ix)a_jb_j + b_j(1 + a_i + x + a_ix) \in I \\
1 + b_i + x + b_ix &= x(1 + a_i + x + a_ix)x^{-1} + b_i(1 + a_i + x + a_ix) \in I \\
1 + b_i + b_j + b_i b_j &= x(1 + a_i + a_j + a_i a_j)x^{-1} + (1 + a_i + b_i + a_i b_i)\\
& \phantom{=.} + (1 + a_i + a_j + a_ia_j)b_ib_j + (1 + a_i + b_i + a_ib_i)b_j \\ & \phantom{=.}  + (1 + a_j + b_i + a_jb_i)b_j + (1 + a_i + b_j + a_i b_j) \in I.
\end{aligned}
\] We have $C_G(A'') = \langle A'', x^2 \rangle$ by Proposition \ref{centralizerofSDP}. From these facts it follows that
\[ I = (1 + s + t + st\, | \, s \in C_G(A''), t \in G).\] Further, every endomorphism of $G$ maps $C_G(A'')$ into itself because $C_G(A'')$ is precisely the subset of elements of order at most 2. By Proposition \ref{centralizercriterion}, $I$ is $\mathrm{End}(G)$-invariant. 

Using the relations in $I$ given above, every element of $R/I$ has the form \[\theta = \alpha + \beta h + \gamma x + \delta x^2\] where $\alpha, \beta, \gamma, \delta \in \{0, 1\}$ and $1 \neq h \in A''.$ Next, we will prove that if $\theta$ is a unit, then $\theta$ is congruent modulo $I$ to an element in  $G$.

Suppose $\theta$ is a unit. If $\beta = 0$, then $\theta^4 = \alpha + \gamma + \delta = 1$. Either each of $\alpha, \gamma, \delta$ is 1 or exactly one is equal to 1. In the former case, $\theta = 1 + x + x^2 = x^3$. In the latter case, $\theta = 1, x$ or $x^2$. Now assume $\beta = 1$. If exactly one of $\alpha, \gamma, \delta$ is zero, then the cases to consider are: \[
\begin{aligned}
1 + a + x &= ax \\
1 + a + x^2 &= ax^2 \\
a + x + x^2 &= a(1 + ax + ax^2) = a(1 + 1 + a + x + 1 + a + x^2) = ax^3.
\end{aligned}
\] The remaining elements to consider are: $1 + a, a + x, a+ x^2, 1 + a + x + x^2$. These elements have augmentation zero and are therefore not units.

Finally, we must check that $1 + u \not \in I$ for $1 \neq u \in G$. Write $u = hx^\epsilon$ where $h \in A''$. The quotient map $G \longrightarrow \C_4$ induces a surjective ring homomorphism \[ R/I \longrightarrow \F_2[\C_4].\] Thus, if $\epsilon \neq 0,$ then $1 + u \not \in I$. So assume $u = h \in A''$. The $Q$-submodule $M$ generated by $u$ is isomorphic to $Q$; since this ring is quasi-frobenius, $M$ is a summand of $A''$ and we may write $A'' = M \times M'$ where $\C_4$ acts decomposably on the direct product. Thus, there is a surjective group homomorphism $G \longrightarrow Q \rtimes \C_4$ that maps $u$ to a nonidentity element. This in turn induces a surjective ring homomorphism \[R/I \longrightarrow \F_2[Q \rtimes \C_4]/(1 + a + x + ax, 1 + x + x^2 + x^3).\] It therefore suffice to prove that $\F_2[Q \rtimes \C_4]/(1 + a + x + ax, 1 + x + x^2 + x^3)$ realizes $Q \rtimes \C_4$ (this group has GAP id (16,3)); this was done using GAP (\cite[\href{https://cocalc.com/klockrid/main/FPENAG/files/Order-16.ipynb}{Order-16.ipynb}]{ourcode}
\end{proof}

\begin{thm} \label{hc4thm}
Let $A$ be a finite elementary abelian $2$-group equipped with a $\C_4= \langle x \rangle$-action. The group $A \rtimes \C_4$ is fully realizable if and only if $x^2$ acts trivially.
\end{thm}
\begin{proof}
The `if' direction follows from Proposition \ref{hc4} and Proposition \ref{sdpxh}: in \ref{sdpxh}, we may ignore condition (\ref{a0}) since we need only apply the proposition when $H = \C_2^L$, condition (\ref{a1}) is obviously satisfied, and condition (\ref{a2}) holds since every element of order 2 lies in $C_G(A'')$. 

For the `only if' direciton, we must show that when the action of $\C_4$ on $A$ is faithful, $A \rtimes \C_4$ is not fully realizable. So, suppose the action of $\C_4$ on $A$ is faithful. Below, we summarize the action of $\C_4$ on $U = S/(z^3)$ and $S$ using multiplicative notation. Let $x = 1 + z$. We have
\[
\begin{aligned}
U &= \C_2^3 = \langle a \rangle \times \langle b \rangle \times \langle c \rangle &&
a \mapsto a^x = ab,
b \mapsto b^x = bc,
c \mapsto c^x = c \\
S &= \C_2^4 = \langle a \rangle \times \langle b \rangle \times \langle c \rangle \times \langle d \rangle &&
a \mapsto a^x = ab,
b \mapsto b^x = bc,
c \mapsto c^x = cd,
d \mapsto d^x = d.
\end{aligned}
\] If $A \rtimes \C_4$ is fully realizable, then either $U \rtimes \C_4$ or $S \rtimes \C_4$ is fully realizable. The first group has GAP id (32, 6), which we determined was realizable but not fully realizable using GAP (\cite[\href{https://cocalc.com/klockrid/main/FPENAG/files/SG-32-6.ipynb}{SG-32-6.ipynb}]{ourcode}). The second group $G$ has GAP id (64, 32). Using GAP, we found an element $a$ in $G$ with $|a| \geq 8$ and $C_G(a) = \langle a \rangle$ (\cite[\href{https://cocalc.com/klockrid/main/FPENAG/files/SG-64-32.ipynb}{SG-64-32.ipynb}]{ourcode}). By \cite[1.4]{SW20}, $G$ is not realizable in characteristic 2.
\end{proof}

The nonabelian, indecomposable fully realizable groups in this subsection may be inductively constructed using the module $Q$: $G_1 = Q \rtimes \C_4$ and \[G_i = \underbrace{Q \rtimes  Q \rtimes \cdots \rtimes Q}_i\rtimes \C_4\] for $i \geq 2$.

\subsection[Case: n = 3]{Cyclic of order 3} \label{sec:n3}

In this section, our goal is slightly more ambitious. We will prove the following theorem.

\begin{thm} \label{na2sg}
Let $G$ be a finite group with abelian normal Sylow $2$-subgroup. The group $G$ is fully realizable if and only if $G = W \times (H \rtimes \C_3),$ where $W$ is a subgroup of $\C_4$ and $H$ is a finite elementary abelian $2$-group.
\end{thm}

\begin{prop} \label{normSyl2}
Let $G$ be a finite group of even order that is fully realizable, and let $H$ be a Sylow $2$-subgroup of $G$. If $H$ is normal in $G$,  then $G$ is either $H$ or $H \rtimes \C_3$. In particular, $|G|$ is either $2^n$ or $3(2^n)$ for some $n \ge 1$.
\end{prop}

\begin{proof}
The normal Sylow $2$-subgroup $H$,  and its index $[G \colon H]$ in $G$, are relatively prime. Therefore, by the Schur-Zassenhaus theorem, we know that $H$ has a complement $L$ in $G$. This means $G = HL$ and $H \cap L = \{ e\}$. In other words, $G = H \rtimes L$.  Moreover, $L$ is a retract of $G$ via the quotient map 
\[ G \xrightarrow{} G/H \cong L.\]
Since $G$ is fully realizable, so is $L$. But since $L$ is a group of odd order, it must be either $\C_1$ or $\C_3$. This completes the proof of this proposition.
\end{proof}

By Proposition \ref{normSyl2}, if the Sylow 2-subgroup $H$ of a fully realizable group $G$ is normal, then either $G = H$ or $G = H\rtimes \C_3$. If $G$ is abelian, then in \cite{chebolu2022fuchs} we proved that $H$ has the form indicated in Theorem \ref{na2sg} (the semidirect product would be a direct product in this case). So, to prove Theorem \ref{na2sg}, it suffices to consider the case where $G = H \rtimes \C_3$. The semidirect product is determined by an action of $\C_3 = \langle c \rangle$ on $H$ via an automorphism $h \mapsto h^c$. By Proposition \ref{fpfsplit}, we may reduce to the case where this action is fixed-point free.

Now we must prove two things: (1) for any elementary abelian $2$-group $H$ and any action of $\C_3$ on $H$, $\C_4 \times (H \rtimes \C_3)$ is fully realizable; and (2) if the action of $\C_3$ on the finite abelian $2$-group $H$ is fixed-point free and $H \rtimes \C_3$ is fully realizable, then $H$ is elementary. Let's begin with task (1).

\begin{prop}
Suppose the action of $\C_3$ on the finite $2$-group $A$ is fixed-point free. Then, for all $a \in Z(A)$ there exists $u  \in Z(A)$ such that $a = u^{c}u^{-1}.$
\end{prop}
\begin{proof}
Define $\psi \colon Z(A) \longrightarrow Z(A)$ by $\psi(u) = u^{c}u^{-1}.$ The map $\psi$ is a group homomorphism since $Z(A)$ is abelian. If $u \in \ker \psi$, then
\[ u^{c}u^{-1} = 1 \implies u^{c} = u  \implies u = 1.\]
Thus, $\psi$ is injective. Since $A$ is finite, $\psi$ is an isomorphism. This completes the proof.
\end{proof}

\begin{prop} \label{aiaj}
Suppose the action of $\C_3$ on the finite $2$-group $A$ is fixed-point free. Then, for all $a_1, a_2 \in A$ with $a_2 \in Z(A)$, \[ 1 + a_1 + a_2 + a_1a_2 \in (1 + a_1 + c + a_1c)\subseteq \F_2[A \rtimes \C_3].\] 
\end{prop}
\begin{proof}
Let $G = A \rtimes \C_3$ and assume the action of $\C_3$ on $A$ is fixed-point free. Note that $ca = a^c c$ for $a \in A$. Take any $a_1, a_2 \in A$  with $a_2 \in Z(A)$ and use the previous proposition to pick $u \in Z(A)$ such that $a_2 = u^{c}u^{-1}$ Then,
\[
\begin{aligned}
u(1 + a_1 + c + a_1c)u^{-1} + (1 + a_1 + c + a_1c)a_2 &= 1 + a_1 + u(u^{-1})^c c+ a_1u(u^{-1})^c c \\ & \phantom{ee} + a_2 + a_1 a_2 + a_2^cc + a_1 a_2^c c \\
&= 1 + a_1 + a_2 + a_1 a_2 \\
& \phantom{ee} + (u(u^{-1})^c + a_1u(u^{-1})^c + a_2^c + a_1a_2^c)c \\
&= 1 + a_1 + a_2 + a_1 a_2 \\
& \phantom{ee} + (a_2^c + a_1a_2^c + a_2^c + a_1a_2^c)c \\
&= 1 + a_1 + a_2 + a_1 a_2.
\end{aligned}
\]
\end{proof}

\begin{prop} \label{whenunit}
Let $G$ be a finite group with normal Sylow $2$-subgroup $H$ of index $3.$ Let \[F \colon \F_2[G] \longrightarrow \F_2[\C_3]\] be the surjective ring homomorphism induced by the quotient map $G \longrightarrow G/H = \C_3$. Then,  $\tau \in \F_2[G]$ is a unit if and only if $F(\tau)$ is a unit.
\end{prop}
\begin{proof}
The `only if' direction follows from the fact that $F$ is a ring homomorphism. 

For the `if' direction, suppose $u = F(\tau)$ is a unit. Since $F$ is surjective, there is an elemnt $\sigma \in \F_2[G]$ such that $F(\sigma) = u^{-1}$. Since $H$ is a $2$-group, the Jacobson radical of $\F_2[H]$ is the augmentation ideal $A$ of $\F_2[H]$. By \cite[16.6]{P71}, the Jacobson radical of $\F_2[G]$ is $J = A\F_2[G]$. The kernel of the map $F$ defined in the statement of the proposition is exactly $J$. Now, \[F(1 + \tau\sigma) = 1 + uu^{-1} = 1 + 1 = 0.\] Thus, $1 + \tau\sigma \in J$ and so $\tau\sigma$ is a unit, implying that $\tau$ is a unit as well (in a finite ring, every element is either a unit or a zero divisor).
\end{proof}

\begin{prop}
Let $A$ be a finite elementary abelian $2$-group with $\F_2$-vector space basis $\{a_1$, $\dots, a_n\}$ and suppose $A$ is equipped with a fixed point-free $\C_3$-action. Let $I$ denote the two-sided ideal of $\F_2[A \rtimes \C_3]$ generated by the elements $1 + a_i + c + a_ic$ for $i = 1, \dots, n$.  The ring $\F_2[A \rtimes \C_3]/I$ fully realizes $A \rtimes \C_3$.\label{sdp}
\end{prop}
\begin{proof}
Let $G = A \rtimes \C_3$. Take $s \in C_G(A) = A$ and $t \in G$. Using Proposition \ref{aiaj} and the definition of $I$, one can show that $1 + s + t + st \in I$. Further, any endomorphism of $G$ maps $A$ into $A$. Hence, by Proposition \ref{centralizercriterion}, the ideal $I$ is $\mathrm{End}(G)$-invariant.

We next prove that $\F_2[G]/I$ realizes $G$. To do so, we will prove that every unit in $\F_2[G]/I$ is congruent to an element of $G$ (modulo $I$) and that $G$ embeds in $\F_2[G]/I$.

Every element of $\F_2[G]/I$ has the form $\theta = \alpha + \beta a + \gamma c + \delta c^2$ where $\alpha, \beta, \gamma, \delta \in \{0,1\}$ and $a \in A$ with $a \neq 1$. Suppose $\theta$ is a unit. We must show that $\theta$ is congruent to an element of $G$ modulo $I$. First, suppose $\beta = 0$. Then, \[
\begin{aligned}
\theta^2 &= (\alpha + \gamma c + \delta c^2)^2 \\
& = \alpha + \gamma c^2 + \delta c \\
\theta^4 &= \alpha + \gamma c + \delta c^2 \\
&= \theta.
\end{aligned}
\] Hence, if $\theta$ is a unit, then $\theta^3 = 1$. Now, \[
\begin{aligned}
1 = \theta^3 &= (\alpha + \gamma c^2 + \delta c)(\alpha + \gamma c + \delta c^2) \\
&= (\alpha + \delta + \gamma) + (\alpha \gamma + \alpha \delta + \gamma \delta)c + (\alpha \delta + \alpha \gamma + \gamma \delta)c^2. 
\end{aligned}\] This forces $\alpha + \delta + \gamma = 1$, meaning either exactly one of $\alpha, \delta, \gamma$ is 1 or all three are equal to 1. In the first case, one obtains $\theta = c^\epsilon$ for some $\epsilon$. The other case yields $1 = 1 + c + c^2$, a contradiction. Now suppose $\beta = 1$. If exactly one of $\alpha, \gamma, \delta$ is zero then $\theta = 1 + a + c = ca$, or $\theta = 1 + a + c^2 = c^2 a$, or $\theta a = (a + c + c^2)a = (1 + c + c^2)$. In the first two cases, $\theta$ is a trivial unit; in the last case $\theta$ is not a unit because $1 + c + c^2$ is not a unit (this was established above, when $\beta = 0$). It remains to consider $1 + a, a + c, a + c^2,$ and $1 + a + c + c^2$. Under the map $F$ in Proposition \ref{whenunit}, these elements map to $0, 1 + c, 1 + c^2,$ and $c + c^2$, none of which are units in $\F_2[\C_3]$, and so none of these elements are units.

Finally, we must check that $G$ embeds in $\F_2[G]/I$; i.e., we must show that $1 + u \not \in I$ for $1 \neq u \in G$. Write $u = ac^\epsilon$ where $a \in A$. The map $F$ from Proposition \ref{whenunit} satisfies $F(I) = 0$, so $F$ descends to a surjective ring homomorphism $\F_2[G]/I \longrightarrow \F_2[\C_3]$. Thus, if $\epsilon \neq 0,$ then $1 + u \not \in I$. So assume $u = a \in A$. Since $A$ is an $\F_2[\C_3] = \F_2 \times \F_4$-module and the action is fixed-point free, $u \in A$ generates a submodule $M$ of $A$ that is isomorphic to $\F_4$. Since $\F_2[\C_3]$ is quasi-frobenius, the projective submodule $M$ is a summand, and we may write $A = M \times M'$, where $\C_3$ acts decomposably on the direct sum. Thus, there is a surjective group homomorphism $G \longrightarrow \F_4 \rtimes \C_3$ that maps $u$ to a non-identity element. This in turn induces a surjective ring homomorphism \[\F_2[G]/I \longrightarrow \F_2[\F_4 \rtimes \C_3]/(1 + u + c + uc, 1 + u^c + c + u^cc).\] It therefore suffices to prove that $\F_2[\F_4 \rtimes \C_3]/(1 + u + c + uc, 1 + u^c + c + u^cc)$ realizes $\F_4 \rtimes \C_3$ (this group is $\A_4$); this was done using GAP (\cite[\href{https://cocalc.com/klockrid/main/FPENAG/files/A4.ipynb}{A4.ipynb}]{ourcode}).
\end{proof}

\begin{prop}
For any finite elementary abelian $2$-group $H$, $G = \C_4 \times (H\rtimes \C_3)$ is fully realizable. 
\end{prop}
\begin{proof}
Note that $G$ has the form $\C_4 \times B \times (A \rtimes \C_3)$ where the action of $\C_3$ on the elementary abelian 2-group $A$ is fixed-point free and $B$ is a finite elementary abelian $2$-group. By Proposition \ref{sdp}, $A \rtimes \C_3$ is fully realizable by $\F_2[A\rtimes \C_3]/J$ where $J = (1 + c + a_i + a_i c)$ and the $a_i$ vary over a basis for $A$. In the proof of that proposition, it is also established that $1 + a + b + ab \in I$ for all $a, b \in A$. Further, note that if $a, b \in A \rtimes \C_3$ are commuting elements with $a \in A$, then $b \in A$ as well. Thus, the conditions in Proposition \ref{sdpxh} are satisfied and $G$ is fully realizable.
\end{proof}

The nonabelian, indecomposable fully realizeable groups in this subsection are determined by the unique 2-dimensional $\F_2[\C_3] = \F_2 \times \F_4$-module on which $\C_3$ acts faithfully. This module is $\F_4$ with module structure determined by the quotient map $\F_2 \times \F_4 \longrightarrow \F_4$. Using multiplicative notation, where $x$ generates $\C_3$, this module $Y$ may be summarized as follows:
\[
\begin{aligned}
Y &= \C_2^2 = \langle a \rangle \times \langle b \rangle && a \mapsto a^x = ab, b\mapsto b^x = a.
\end{aligned}
\]
The mentioned examples may be constructed inductively: $G_1 = Y \rtimes \C_3 = \A_4$ and \[G_i = \underbrace{Y \rtimes  Y \rtimes \cdots \rtimes Y}_i\rtimes \C_3\] for $i \geq 2$

Now we turn to our last task: show that if the action of $\C_3$ on a finite abelian 2-group $H$ is fixed-point free and $H \rtimes \C_3$ is fully realizable, then $H$ is elementary.

\begin{prop}
Let $H$ be a finite $2$-group on which $\C_3 = \langle c \rangle$ acts via an automorphism $h \mapsto h^c$ whose restriction to $Z(H)$ is fixed-point free. If $H\rtimes \C_3$ is realizable, then $Z(H)$ is elementary abelian. In particular, if $H$ is abelian then it is elementary abelian.
\end{prop}
\begin{proof}
Let $G = H \rtimes \C_3$ and take $h \in Z(H)$. It suffices to prove that $|h| = 1$ or $2$. By Proposition \ref{whenunit}, $1 + h + c$ is a unit in $\F_2[G]$. Suppose $R = \F_2[G]/I$ realizes $G$. In $R$, we must have $1 + h + c = t \in G$. If $t \in H$, then $c$ lies in the subring $S$ of $R$ generated by $H$, which forces the order of $c$ to be a power of 2, a contradiction. Hence, $t = ca$ or $t = c^2a$ for some $a \in H$.

Suppose $1 + h + c + ca = 0$ in $R$. Then \[\theta = c^2(1+h) = c^2(c(1 + a)) = c^3(1+a) = 1 + a \in S,\] and hence $\theta$ commutes with $h$. Thus,
\[
\begin{aligned}
hc^2(1+h) &= c^2(1+h)h \\
c^2 h^c(1 + h) &= c^2(1+h)h \\
h^c(1 + h) &= (1+h)h \\
h^c + h^ch &= h + h^2 \\
1 + h + h^c + h^ch^{-1} &= 0.
\end{aligned}
\] The element $\tau = h + h^c + h^{c^2} \in \F_2[G]$ is a central unit with $\gcd(|\tau|, 3) = 1$ and $\tau^c = \tau$. In $R$, we must therefore have that $\tau$ is congruent to a unit in $Z(H)$ that is fixed by the $\C_3$ action. This forces $\tau = 1$:\[ 1 + h + h^c + h^{c^2} = 0.\] Further, since the action of $\C_3$ on $Z(H)$ is fixed-point free, $hh^ch^{c^2} = 1$. Combining this with the work above, we obtain
\[
\begin{aligned}
h^ch^{-1} &= h^{c^2} \\
h^c &= hh^{c^2} \\
h^c &= (h^c)^{-1} \\
(h^c)^2 &= 1.
\end{aligned}
\] This implies that $|h| = |h^c| = 1$ or 2. The case where $t = c^2a$ is similar.
\end{proof}

\subsection{Nonabelian 2-groups: a family of examples} \label{sec:nab2}

In the remainder of this section, we construct one more family of fully realizable groups in order to show that there are examples of nonabelian $2$-groups $H$ where $H \rtimes \C_3$ is fully realizable. To do so, we will consider $\C_6$-actions on elementary abelian $2$-groups $A$; note that $A \rtimes \C_6 = (A\rtimes \C_2)\rtimes \C_3$.
Above, we defined indecomposable $\F_2[\C_2]$-modules:
\[
\begin{aligned}
\F_2 &= \C_2 = \langle a \rangle && a \mapsto a\\
Q &= \C_2^2 = \langle a, b \rangle && a \mapsto ab, b \mapsto b \\
\end{aligned}
\]
and indecomposable $\F_2[\C_3]$-modules:
\[
\begin{aligned}
\F_2 &= \C_2 = \langle a \rangle && a \mapsto a\\
Y &=\F_4 = \C_2^2 = \langle a, b \rangle && a \mapsto ab, b \mapsto a. \\
\end{aligned}
\] Now, $\F_2[\C_6] = \F_2[\C_2] \otimes \F_2[\C_3] = Q \times Y \otimes Q$. Every finitely generated $\F_2[\C_6]$-module is a sum of the following modules:
\[
\begin{aligned}
\F_2 &= \C_2 && \text{$\C_2, \C_3$ act trivially}\\
Q&= \C_2^2 &&\text{faithful $\C_2$-action, $\C_3$ acts trivially}\\
Y &= \C_2^2 &&\text{faithful $\C_3$-action, $\C_2$ acts trivially} \\
Y \otimes Q &= \C_2^4 = \langle \alpha, \beta, \gamma, \delta \rangle && \text{tensor product of faithful $\C_3$- and $\C_2$-actions} \\
& && \alpha \mapsto \alpha\beta\gamma\delta, \beta \mapsto \beta \delta, \gamma \mapsto \alpha\beta, \delta \mapsto \beta \,\, (\C_6\text{-action})\\
& && \alpha \leftrightarrow a \otimes a, \beta \leftrightarrow a \otimes b, \gamma \leftrightarrow b \otimes a, \delta \leftrightarrow b \otimes b.
\end{aligned}
\]

\begin{thm}
Let $A$ be a finite elementary abelian $2$-group equipped with a $\C_6 = \langle x \rangle$-action. The group $A \rtimes \C_6$ is fully realizable if and only if the module $Q$ does not appear as a summand of $A$. \label{c6}
\end{thm}

\begin{proof} Let $G = A \rtimes \C_6$. To prove the `only if' direction, we may decompose $A$ as a direct sum of the modules $\F_2, Q, Y$ and $Y \otimes Q$ defined above. If $A$ contains $Q$ as a summand and $A \rtimes \C_6$ is fully realizable, then so is $Q \rtimes \C_6 = \C_3 \times (Q \rtimes \C_2)$, but as mentioned above $Q \rtimes \C_2$ is the dihedral group of order 8 which is not fully realizable.

For the `if' direction, write $A = \F_2^{m''} \times Y^{m'} \times (Y \otimes Q)^m$. Since $\C_6$ acts decomposably on this direct sum decomposition of $A$ and trivially on $\F_2^{m''}$, we may assume $m'' = 0$ and then use Proposition \ref{sdpxh} to recover this term. Since $m'' = 0$, the $\C_6$-action on $A$ is fixed-point free.

Write $\alpha_i, \beta_i, \gamma_i, \delta_i$ for the generators of the $i$th copy of $Y \otimes Q$ in $A$. Define two subgroups $A_1, A_2 \leq A$:
\[
\begin{aligned}
A_1 &= Y^{m'} \times \oplus_{i=1}^m \langle \beta_i, \delta_i\rangle \\ 
A_2&= \oplus_{i=1}^m  \langle \alpha_i, \gamma_i\rangle.
\end{aligned}
\] Note that $A = A_1 \times A_2$. One can check that the following identities hold for the generators of each summand $Y \otimes Q$:
\[
\begin{aligned}
\alpha^{x^3} &= \alpha\beta &\beta^{x^3} &= \beta &\gamma^{x^3} &= \gamma \delta& \delta^{x^3} &= \delta \\
\alpha^{x^2} &= \gamma  &\beta^{x^2} &= \delta &\gamma^{x^2} &= \alpha\gamma  &\delta^{x^2} &= \beta\delta \\
\end{aligned}
\] The $\C_6$-action maps $A_1$ to itself and $x^3$ acts trivially on $A_1$. The action of $x^2$ maps $A_2$ to itself. Since \[ (ax^\epsilon)^k = a a^{x^\epsilon}\cdots a^{x^{\epsilon(k-1)}} x^{\epsilon k},\] we have that $|x^\epsilon|$ divides $|ax^\epsilon|$. In particular, if $\epsilon = 1, 3$ or $5$, then $ax^\epsilon$ has even order, so an element of $G$ of order 3 has the form $ax^{2\epsilon}$ where $\epsilon \neq 0$. Since $a a^{x^2} a^{x^4} = 1$ for all $a \in A$ (one can directly check that the identity holds for each generator of $A$), $ax^{2\epsilon}$ has order 3 for any $a \in A$.

Let $V = A \rtimes \langle x^3 \rangle$ and let $W = A \rtimes \langle x^2 \rangle$. We have \[ G = V \rtimes \langle x^2\rangle = W \rtimes \langle x^3 \rangle.\] For any endomorphism $\rho$ of $G$, we have $\rho(V) \subseteq V$ since $V$ is a 2-group of index 3 in $G$. Further, $\rho(A) \subseteq A$: assume to the contrary that there exists an $a \in A$ such that $\rho(a) = a'x^3$. Since $x^2$ has order 3, we may write $\rho(x^2) = a''x^{2\epsilon}$. Now, 
\[
\begin{aligned}
\rho(ax^2) &= a' x^3 a'' x^{2\epsilon} \\
&= a' (a'')^{x^3} x^{3 + 2\epsilon}.
\end{aligned}
\] Above, $2 + 3\epsilon = 1, 3,$ or $5$, implying that $|\rho(ax^2)|$ is 2, 4, or 6. But $|\rho(ax^2)| = 3$ since $ax^2$ has order 3, a contradiction.\footnote{If $\rho$ were just an endomorphism of $V$, it would not necessarily take $A$ into itself!} Since $\rho(x^2)$ has order 1 or 3 and $\rho(A) \subseteq A$, we have $\rho(W) \subseteq W$.

Let $I\subseteq \F_2[G]$ be the ideal generated by elements of the form \[ 1 + v + wx^{3\epsilon} + v^{x^{3\epsilon}} wx^{3\epsilon} \] where $\epsilon \in \{0, 1\}$, $v \in V$ and $w \in W$. We claim that $I$ fully realizes $G$.

First, we prove that $I$ is $\mathrm{End}(G)$-invariant. Take $\rho \in \mathrm{End}(G)$. For $\epsilon = 0$, \[ 1 + v + w + vw \mapsto 1 + \rho(v) + \rho(w) + \rho(v)\rho(w) \in I\] since $\rho(V) \subseteq V$ and $\rho(W) \subseteq(W)$. Now suppose $\epsilon = 1$ and consider two cases. Since $x^3$ has order 2, we have $\rho(x^3) = a_1 x^3$ or $\rho(x^3) = a$, where $a \in A$ and $a_1 \in A_1$. In the former case:
\[
\begin{aligned}
1 + v + wx^3 + x^3vx^3wx^3 &\mapsto 1 + \rho(v) + \rho(w)a_1x^3 + a_1x^3\rho(v)a_1x^3\rho(w)a_1x^3 \\
&= 1 + a_1 \rho(v) a_1 + \rho(w)a_1x^3 + (a_1 \rho(v) a_1)^{x^3} \rho(w)a_1x^3 \in I
\end{aligned}
\] (we used that $a_1 \in Z(V)$). In the latter case, note that $\rho(x^3) =a \in A$ implies $\rho(v) \in A$, so $a$ commutes with $\rho(v)$. Thus, 
\[
\begin{aligned}
1 + v + wx^3 + x^3vx^3wx^3 &\mapsto 1 + \rho(v) + \rho(w)a + a\rho(v)a\rho(w)a \\
&= 1 + \rho(v) + \rho(w)a + \rho(v)\rho(w)a \in I.
\end{aligned}
\] This completes the proof that $I$ is $\mathrm{End}(G)$-invariant.

Next, we prove that every unit in $\F_2[G]/I$ is congruent to an element of $G$. Note that the following equations hold modulo $I$:
\[
\begin{aligned}
ax^\epsilon &= 1 + a^{x^{3\epsilon}} + x^\epsilon && \text{ for $a \in A$ and $0 \leq \epsilon \leq 5$} \\
x^{3+k} &= 1 + x^3 + x^k && \text{ for $k = 1, 2$.}
\end{aligned}
\] Thus, every element of $\F_2[G]/I$ has the form \[ \theta = c_0 + c_1a + c_2 x + c_3x^2 + c_4x^3\] where $a \in A$ and each $c_i \in \{0, 1\}$. By Proposition \ref{whenunit}, $\theta$ is a unit if and only if exactly one of $c_0 + c_1 + c_4$, $c_2$, $c_3$ is 1 and the rest are 0.
If $c_0 + c_1 + c_4 = 1$ and $c_2 = c_3 = 0$, then the possible values of $\theta$ are $1, a, x^3, 1 + a + x^3 = x^3a$. If $c_0 + c_1 + c_4 = c_2 = 0$ and $c_3 = 1$, then the possible values of $\theta$ are $x^2, x^2 + a + x^3 = x^5, x^2 + 1 + a = ax^2, x^2 + a + x^3 = x^3ax^3.$ If $c_0 + c_1 + c_4 = c_3 = 0$ and $c_2 = 1$, then the possible values of $\theta$ are $x, x + 1 + x^3 = x^4, x + 1 + a = a^{x^3}x, x + a + x^3 = x^3 a^{x^3}x.$ In all cases we obtain that $\theta$ is a trivial unit.

Finally, we must check that no element of $G$ is congruent to the identity of $G$ modulo $I$. To do so, it suffices to explicitly prove full realizability for the following values of $A$: $A = \F_2, Y, Y \otimes Q$. The corresponding groups are $G = \C_2 \times \C_6, \C_2 \times \A_4$, and $\C_2^4 \rtimes \C_6$, where the last group has GAP id (96, 70). In each of these cases, we used GAP to check that the ideal $I$ defined above indeed realizes $G$ (\cite[\href{https://cocalc.com/klockrid/main/FPENAG/files/C6-Cases.ipynb}{C6-Cases.ipynb}]{ourcode}).

To recover the term $\F_2^{m''}$, we must verify the conditions in Proposition \ref{sdpxh}. The first condition does not apply, and the second condition is clearly satisfied. For condition (\ref{a2}), note that if $v$ has order 2 then $v = a \in A$ or $v = a_1x^3$ where $a_1 \in A_1$. In the latter case, $x^3$ fixes $v$ and $1 + v + g + vg \in I$ for any $g \in G$. In the former case, take $g = wx^{3\epsilon} \in G$ where $w = a'x^{2\tau}\in W$ and suppose $v$ commutes with $g$. If $\epsilon = 0$ then $1 + v + g + vg \in I$. If $\epsilon = 1$, then $a = a^{x^{2\tau+ 3}}$. If $\tau = 0$ then $x^3$ fixes $v$ and again $1 + v + g + vg \in I$. If $\tau = 2$ or $4$, then $a = 1$. This completes the proof.
\end{proof}


\bibliographystyle{alpha}
\bibliography{fnab.bib}

\end{document}